\documentclass[11pt]{amsart}
\usepackage[margin=1in]{geometry}
\usepackage{bbm}
\usepackage{float, graphicx}
\usepackage[]{epsfig}
\usepackage{amsmath, amsthm, amssymb,}
\usepackage{epsfig}
\usepackage{verbatim}
\usepackage{multicol}
\usepackage{url}
\usepackage{latexsym}
\usepackage{mathrsfs}
\usepackage[colorlinks, bookmarks=true]{hyperref}
\usepackage{graphicx}
\usepackage{amsmath}
\usepackage{enumerate}
\usepackage[normalem]{ulem}
\usepackage{bm}
\usepackage{dsfont}
\usepackage{stmaryrd}
\usepackage{enumitem}
\usepackage{soul}
\usepackage[noadjust]{cite}
\usepackage{color}
\usepackage{xcolor}
\usepackage{hyperref}
\usepackage{indentfirst}
\usepackage[capitalise]{cleveref}

\usepackage{mathtools}

\newcommand{\cG}{{\mathcal G}}

\newcommand{\bR}{{\mathbb R}}

\newcommand{\bP}{{\mathbb P}}

\def\<{\langle}
\def\>{\rangle}

\def\cG{\mathcal{G}}

\def\Lamp[#1]{\boldsymbol{\Lambda}_{\mathrm{AMP}}^{(#1)}}
\def\lalg[#1]{\Lambda_{\mathrm{alg}, #1}}

\def\p{\mathbb{P}}

\newcommand{\oo}{{\rm o}}

\newcommand{\cor}{\color{red}}

\newcommand{\Tr}{{\rm Tr}}

\newcommand{\E}{\mathbb{E}}

\newcommand{\Var}{\operatorname{Var}}

\newcommand{\norm}[1]{\left\lVert#1\right\rVert}
\newcommand{\RN}[1]{%
  \textup{\uppercase\expandafter{\romannumeral#1}}%
}

\newcommand{\RNum}[1]{\uppercase\expandafter{\romannumeral #1\relax}}

\theoremstyle{plain} %plain, definition, remark
\newtheorem{theorem}{Theorem}[section]
\newtheorem*{theorem*}{Theorem}
\newtheorem{lemma}[theorem]{Lemma}
\newtheorem*{lemma*}{Lemma}

\newtheorem*{corollary*}{Corollary}
\newtheorem{proposition}[theorem]{Proposition}
\newtheorem*{proposition*}{Proposition}

\newtheorem*{assumption*}{Assumption}

\newtheorem{definition}[theorem]{Definition}
\newtheorem*{definition*}{Definition}

\newtheorem*{example*}{Example}
\newtheorem{remark}[theorem]{Remark}

\newtheorem*{remark*}{Remark}
\newtheorem*{remarks*}{Remarks}

\title{The spectral distribution of random graphs with given degree sequences}
\date{}

\begin{document}

\author{Shuyi Wang\textsuperscript{1}}
\address{\textsuperscript{1}University of Pennsylvania, Philadelphia, PA, USA}
\email{sywang25@sas.upenn.edu}

\author{Kevin Li\textsuperscript{2}}
\address{\textsuperscript{1}University of Pennsylvania, Philadelphia, PA, USA}
\email{kevinmli@seas.upenn.edu}    

\author{Jiaoyang Huang\textsuperscript{3}}
\address{\textsuperscript{3}University of Pennsylvania, Philadelphia, PA, USA}
\email{huangjy@wharton.upenn.edu}

\begin{abstract}
In this article, we study random graphs with a given degree sequence $d_1, d_2, \cdots, d_n$ from the configuration model. We show that under mild assumptions of the degree sequence, the spectral distribution of the normalized Laplacian matrix of such random graph converges in distribution to the semicircle distribution as the number of vertices $n\rightarrow \infty$. This extends work by McKay \cite{Mc} and Tran, Vu and Wang \cite{TVW} which studied random regular graphs ($d_1=d_2=\cdots=d_n=d$). Furthermore, we extend the assumption to show that a slightly more general condition is equivalent to the weak convergence to semicircle distribution. The equivalence is also illustrated by numerical simulations.
\end{abstract}

\maketitle
\section{Introduction} 

\subsection{Backgrounds}

% In random matrix theory, a frequently studied property is the spectral distribution. For random graphs, we can study the spectral distribution of the normalized Laplacian and adjacency matrices. Previous work in this area has mostly focused on settings where the entries (degrees) in the random graph matrix are uncorrelated \cite{CLV} or have short-range correlation \cite{AEK, Ch16}. When given degree sequences are introduced, entries have long-range negative correlation. McKay \cite{Mc} and Tran, Vu and Wang \cite{TVW} studied fixed $d$-regular graphs, showing that the spectral distribution converges to the McKay-Kesten distribution when $d$ is fixed, and to the semicircle distribution when $d\rightarrow\infty$. In this paper, we extend these results to situations where each degree in the degree sequence is allowed to vary, to consider non-regular graphs. In this setting, we show that the spectral distribution of the normalized Laplacian converges in distribution to the semicircle distribution.

Random graphs, since their introduction in the 1960s, have been a widely researched subject that continually receives attentions of scholars. The random graph models shed light on a variety of real-world phenomena and have found applications in the studies of social networks, economics, and biological engineering. Researchers have developed a comprehensive framework of various random graph models, such as the Erd{\H o}s-R{\'e}nyi model, the Preferential Attachment model, the stochastic block model, etc. \cite{bollobas1998random} and \cite{van2024random} are excellent manuscripts that introduce popular random graph models in a systematic manner.

A key tool for studying random graphs is random matrix theory, which often focuses on properties such as spectral distribution. In the context of random graphs, researchers study the spectral distribution of the adjacency matrices and normalized Laplacians, for they provide significant insights into the graph's structural properties.

Previous work in this area has mostly focused on settings where the entries in the random matrix, which are essentially degrees in the random graphs, are uncorrelated \cite{CLV} or have short-range correlation \cite{AEK, Ch16}. When given degree sequences are introduced, entries have long-range negative correlation, which causes difficulties in the study of the spectral distribution.

A plethora of work have been done in the context of random $d$-regular graphs. In \cite{Mc}, Mckay showed that the spectral distribution converges to the Kesten-McKay distribution when $d$ is fixed. \cite{bauerschmidt2019local} showed that the spectral distribution holds under the optimally small scale. \cite{huang2024spectrum} improved on the rigidity of eigenvalues as well as an error bound on the extremal eigenvalue. In \cite{TVW}, Tran, Vu, and Wang showed that the spectral distribution converges to the semicircle distribution when $d\rightarrow\infty$.

In this work, we extend the $d$-regular assumptions to situations where each degree in the degree sequence is allowed to vary, to consider non-regular graphs. Such random graphs are called the configuration model, firstly introduced by Bollob{\'a}s \cite{bollobas1980probabilistic}. The configuration model is frequently studied due to its wide application in statistical inferences and unique properties. In \cite{fosdick2018configuring}, the authors studied the use of configuration model as null model to test empirical graph properties. In \cite{dhara2021critical}, the authors proved a new universality class for the configuration model due to their behavior of critical percolation. More general notions of configuration models are also discussed and investigated in the work of \cite{britton2006generating}, such as erased configuration model and repeated configuration model.

In this paper, we assume mild conditions on the asymptotic behavior of the degree sequence and show that the spectral distribution of the normalized Laplacian converges in distribution to the semicircle distribution. To prove convergence in distribution, we mainly use the moments method. To bound the moments, we use tools from Fourier analysis developed by Sarid \cite{Sa}. To deal with combinatorial difficulties when bounding the moments, we apply several graph transformation techniques. To extend the condition to the equivalent condition, we use birth-death processes to model the graph pruning process. 

\subsection{Preliminaries} \label{sec:preliminaries}
In this section, we introduce the configuration model and establish the necessary notations. Consider a graph with a fixed degree sequence $d_1, d_2, \cdots, d_n$, where the total degree is $D = d_1 + d_2+\cdots + d_n$. The configuration model is constructed as follows: for $n$ vertices $v_1, v_2, \cdots, v_n$, each vertex $v_i$ is associated with $d_i$ half-edges. A random perfect matching is then performed on these $D$ half-edges. Specifically, one half-edge is randomly connected to any of the remaining $D - 1$ half-edges, and the process is repeated for each unused half-edge, reducing the pool of available half-edges by two at each step. This procedure continues until all half-edges are matched. The result is a uniformly random graph (potentially with multi-edges and self-loops) among all graphs with the given degree sequence.

Given a degree sequence $d_1, d_2, \cdots, d_n$, we can sample a graph $\cG$ from the configuration model. This random graph defines an adjacency matrix $A$, where $A_{ij}$ represents the number of edges between vertices $v_i$ and $v_j$. To obtain the normalized Laplacian, we normalize all entries of $A$ to have a mean of 0 and variance $1/n$. Specifically, for each entry $A_{ij}$, we consider each pair of half-edges, one from $v_i$ and one from $v_j$, as an indicator random variable that takes the value 1 if the two half-edges are connected in the random graph. 

More precisely, for each vertex $v_i$, we enumerate its half-edges from 1 to $d_i$. For vertices $v_i$ and $v_j$, and for $1 \leq s \leq d_i$ and $1 \leq t \leq d_j$, we define the indicator function $H_{ij}^{st}$, which equals 1 if the $s$-th half-edge of $v_i$ is connected to the $t$-th half-edge of $v_j$. Then, we can express $A_{ij}$ as
\begin{align}\label{e:indicator_function}
    A_{ij}=\sum_{s=1}^{d_i}\sum_{t=1}^{d_j}H_{ij}^{st}.
\end{align}
 Using linearity of expectation, the mean and variance can be computed as 
 \begin{align}\label{e:A_mean_var}
    \E \left[A_{ij}\right] = \frac{d_id_j}{D-1},\quad \Var \left(A_{ij}\right) \approx \frac{d_id_j}{D}. 
 \end{align}

We define the diagonal matrix $\Delta$, where $\Delta_{ii} = d_i$ represents the degree of vertex $v_i$ in the graph $\mathcal{G}$, and $\Delta_{ij} = 0$ for $i \neq j$. The normalized Laplacian matrix associated with $\mathcal{G}$ is a scaled version of its adjacency matrix $A$, given by:
\begin{align}\label{e:defM}
    M=\sqrt{\frac{D}{n}}\Delta^{-\frac{1}{2}}\left(A-\left[\frac{d_id_j}{D-1}\right]_{ij}\right)\Delta^{-\frac{1}{2}}
    =\sqrt{\frac{D}{n}}\left(\Delta^{-\frac{1}{2}}A\Delta^{-\frac{1}{2}}-\left[\frac{\sqrt{d_id_j}}{D-1}\right]_{ij}\right)
\end{align}
Note that $M$ is a normalized version of the standard Laplacian matrix associated with $\mathcal{G}$ (see \cite{chung1997spectral}). This normalization ensures that $\E[M_{ij}] = 0$ and $\Var(M_{ij}) \approx 1/n$, matching the mean and variance of Wigner matrices. However, unlike Wigner matrices, where the entries are independent, the entries of $M$ exhibit long-range correlations.

Before presenting our main results, we introduce some necessary notation. Let $\mu_n$ denote the empirical eigenvalue distribution of $M$, defined as $\mu_n = \frac{1}{n} \sum_{i=1}^n \delta_{\lambda_i}$, where $\lambda_i$ are the eigenvalues of $M$. We denote by $\rho_{\rm sc}$ the semicircle distribution supported on $[-2, 2]$, with the density function
$$d\rho_{\rm sc} (x) = \frac{1}{2\pi} \sqrt{4 - x^2}dx.$$
The moments of the semicircle distribution are given by the Catalan numbers. Specifically,
$$\frac{1}{2\pi}\int_{-2}^2 x^n \sqrt{4-x^2}dx = \begin{cases}
    0, & \text{    if $n$ is odd,}\\
    C_{n/2}, & \text{    if $n $ is even,}
\end{cases}$$
where $C_k$ is the Catalan number with the following explicit formula
$$C_k = \frac{{2k\choose k}}{k+1}.$$
An important observation is that $C_k$ also has a combinatorial interpretation: it counts the number of distinct rooted trees with $k$ edges.

Throughout the paper, we will use several asymptotic notations. Suppose $\{ a_n\}$ and $\{ b_n\}$ are two sequences of non-negative numbers in $\mathbb{R}$. we write $a_n = o(b_n)$ or $a_n\ll b_n$ if $\lim_{n\rightarrow\infty} {a_n}/{b_n} = 0$;  $a_n = \omega (b_n)$ if $\lim_{n\rightarrow\infty} {a_n}/{b_n} = \infty$; $a_n = O (b_n)$ if there exists a large constant $C>0$, such that $a_n \leq  C b_n$; $a_n = \Omega (b_n)$ if there exists a small constant $c>0$ such that $a_n \geq c b_n$.
 Additionally, for any integer $n$, we use $[n]$ to denote the index set $\{ 1, 2, \cdots , n \}$.

\subsection{Main Results}
In this section, we present the main results of this paper. The first main theorem indicates that moments of the eigenvalue distribution of $M^{(n)}$ converge to those of the semicircle distribution.

\begin{theorem}
    \label{semicircle}
    Consider the degree sequence $\{ d_1^{(n)}, d_2^{(n)}, \cdots, d_n^{(n)} \}$ with total degree $D^{(n)}=d_1^{(n)}+d_2^{(n)}+\cdots+d_n^{(n)}$ for each $n\geq 1$. Suppose that $\min\{d_i^{(n)}\}_{i\in [n]} \gg \sqrt{D^{(n)}/n}$ as $n \rightarrow \infty$. Consider the normalized Laplacian matrix $M^{(n)}$ as defined in \eqref{e:defM}. Then, the moments of the empirical eigenvalue distribution of $M^{(n)}$ converge to the moments of the semicircle distribution:
    \begin{align}
        \label{e:moment_converge}
    \lim_{n\rightarrow \infty}\frac{1}{n}\E[\Tr[({M^{(n)}})^k]] =  
        \begin{cases} 
            C_{k/2} & \text{for even $k$}, \\
            0 & \text{for odd $k$},
        \end{cases}
   \end{align}
    where $C_{k/2}$ is the $k/2$-th Catalan number. Moreover, the empirical eigenvalue distribution $\mu_n$ of $M^{(n)}$ weakly converges to the semicircle distribution $\rho_{\rm sc}$ in probability:
    \begin{align}\label{e:weak_converge}
        \mu_n \xrightarrow{D} \rho_{\rm sc},\quad  \text{in probability}.
    \end{align}
\end{theorem}

\begin{comment}
Therefore, given the two theorems above, we conclude the following weak convergence whose proof follows immediately. 

\begin{theorem}
    \label{weakconv1}
    Consider the degree sequence $\{ d_1^{(n)}, d_2^{(n)}, \cdots, d_n^{(n)} \}$ for each $n\geq 1$. Suppose that $\min\{d_i^{(n)}\}_{i\in [n]} \gg \sqrt{\frac{D^{(n)}}{n}}$ as $n \rightarrow \infty$. The empirical distribution of eigenvalues $\mu_n$ of $M^{(n)}$ weakly converges to the semicircle distribution $\rho_{\rm sc}$ in probability. Notationally,
    $$\mu_n \xrightarrow{D} \rho_{\rm sc} \ \text{in probability}.$$
\end{theorem}
\end{comment}

\begin{remark} \label{rmk:rank1}
   In our statements and proofs, we work with the normalized Laplacian matrix $M$ as defined in \eqref{e:defM},
    $$M=\sqrt{\frac{D}{n}}\Delta^{-\frac{1}{2}}\left(A-\left[\frac{d_id_j}{D-1}\right]_{ij}\right)\Delta^{-\frac{1}{2}}$$
    as it offers intuitive insights and simplifies computations. However, it is also valid to omit the subtraction of the matrix $\left[ \frac{d_i d_j}{D - 1} \right]_{ij}$, since it is rank-$1$ and does not affect the limiting eigenvalue distribution of the resulting matrix $M$.
\end{remark}

In the proof of Theorem \ref{semicircle}, we use the method of moments to establish weak convergence. It is important to note that moment convergence is a sufficient condition for weak convergence. With more refined analysis, the condition ``$\min\{d_i^{(n)}\}_{i \in [n]} \gg \sqrt{D^{(n)}/n}$" in Theorem \ref{semicircle} can be relaxed if we are interested in a more general sufficient condition for weak convergence to the semicircle distribution.

Before introducing this more general condition, we first present a lemma that simplifies the analysis.
\begin{lemma}
    \label{equivlemma}
    We assume that for arbitrarily large $C > 0$ and arbitrarily small $\varepsilon > 0$, the degree sequence satisfies that $\sum_{i = 1}^n \mathbbm{1}\{ d_i^{(n)} < C \sqrt{D^{(n)}/n} \} \leq \varepsilon n$ for all large $n > N(C, \varepsilon)$. Then, with probability $1 - o(1)$, we may remove at most $2 \varepsilon n$ vertices and edges adjacent to them, so that the following holds:
    \begin{enumerate}
        \item The total number of edges removed is no larger than $4\varepsilon C \sqrt{nD}$;
        \item Each vertex of the remaining graph has a degree larger than $C \sqrt{D^{(n)}/n}$;
        \item Condition on the degree sequence of the remaining graph, its law is given by the configuration model.
    \end{enumerate} 
\end{lemma}

Now, we are ready to present the more general sufficient condition for the weak convergence, whose proof follows from Lemma \ref{equivlemma} and \Cref{semicircle}.

% TODO: arbitrary C , eps ; change notation superscript (n)
\begin{theorem} 
% [more general sufficient condition for the weak convergence to semicircle law]
    \label{equiv}
     Consider the degree sequence $\{ d_1^{(n)}, d_2^{(n)}, \cdots, d_n^{(n)} \}$ with total degree $D^{(n)}=d_1^{(n)}+d_2^{(n)}+\cdots+d_n^{(n)}$ for each $n\geq 1$, and recall the normalized Laplacian matrix $M^{(n)}$ as defined in \ref{sec:preliminaries}. Suppose the following condition holds: for arbitrarily large $C > 0$ and arbitrarily small $\varepsilon > 0$, the degree sequence satisfies that $\sum_{i = 1}^n \mathbbm{1}\{ d_i^{(n)} < C \sqrt{D^{(n)}/n} \} \leq \varepsilon n$ for all large $n > N(C, \varepsilon)$. Then the empirical eigenvalue distribution $\mu_n$ of $M^{(n)}$ weakly converges to the semicircle distribution $\rho_{\rm sc}$ in probability: 
    $$\mu_n \xrightarrow{D} \rho_{\rm sc} \quad \text{in probability.}$$
\end{theorem}

\begin{remark}[Necessity of the condition]
    The condition in \Cref{equiv} is also necessary for weak convergence to hold. Suppose there is a fixed proportion of vertices with degrees on the order of $C\sqrt{D^{(n)}/n}$, and this proportion does not vanish as $n \to \infty$. In this case, Lemma \ref{equivlemma} cannot be applied. Furthermore, as noted in Remark \ref{rmk}, we can construct a $k$-step walk that violates moment convergence, meaning that the moments of the eigenvalue distribution will no longer match those of the semicircle distribution. For readability, the technical details of this argument are postponed until Remark \ref{rmk}, after the proof of Theorem \ref{semicircle}.
\end{remark}

From now on, for the simplicity of notations in later proofs, we will drop the superscript of $(n)$ in all symbols.

\subsection{Proof Ideas and Difficulties}
The main difficulty lies in the proof of Theorem \ref{semicircle}. We use a variety of techniques such as the Fourier transform, graph pruning, and combinatorial arguments to tackle the proof. Particularly, the main idea behind is that we expand the quantity in terms of different graph types
$$\frac{1}{n} \E[\Tr M^k] = \sum_G  \frac{1}{n} t_G$$
where each $G$ is a valid closed directed graph with $k$ edges. We basically show that only the graphs that are trees and have $k/2$ edges (thus each edge has multiplicity of $2$) will make a contribution of $1$. Other possible graphs will have contributions of $0$ asymptotically. Therefore, if $k$ is odd, $G$ cannot be a tree and thus this quantity is $0$. If $k$ is even, this sum is counting the number of distinct trees which have $k/2$ edges. That is exactly the combinatorial meaning of the Catalan number $C_{k/2}$.

Although the main idea is straightforward, there are two technical difficulties to overcome in the proof:
\begin{itemize}
    \item[1. ] Fine enough bounds on different graphs are needed so that we can prove the asymptotic contribution of most graphs to be $0$.
    \item[2. ] There are too many possibilities of graphs with $k$ edges. We need to transform most of them into a limited set of graph types and apply combinatorial arguments so that we only need give bounds on a small number of graph types.
\end{itemize}

To tackle the first point, we have to attain a good upper bound on the part of graph where edges have multiplicities of $1$. It is because this part appears most frequently in our computations. To do it, we recreate the methods developed in \cite{Sa}. We develop a Fourier transform on the conditional probabilities of the edge indicators and bound the Fourier coefficients. 

To tackle the second point, we introduce a variety of graph transformations which only increase $t_G$. Therefore, whenever we show that the transformed graph has contributions of $o(1)$, we know that the untransformed graphs have contributions of $o(1)$. We basically define all possible transformations and show that the transformed graph consists of several components whose contributions are easier to compute.

Another difficulty lies in proving Lemma \ref{equivlemma}, since we need to make sure that the pruning operation will not continue forever and lead to an empty set of vertices in the end. We use a birth and death chain to keep track of the number of half-edges in vertices whose degrees are less than the threshold. By iterative probability arguments, we bound the stopping time of the chain and show that the impact will be no more than $o(n)$ vertices.

% \subsection{Organization of the Paper}

% In Section \ref{sec:semicircle}, we present the proof of Theorem \ref{semicircle}.

% the moments convergence of the empirical eigenvalue distribution.

% In Section \ref{sec:concentration}, we present the full proof of Theorem \ref{concentration}, which shows that the empirical eigenvalue distribution will converge in probability.

% In Section \ref{sec:equivlemma}, we present the full proof of Lemma \ref{equivlemma} and \Cref{equiv}. 

% In Section \ref{numerical}, we present numerical simulations to support our findings. 

\section{Proof of Theorem \ref{semicircle}} \label{sec:semicircle}

% In this section we prove the semicircle law for the normalized adjacency matrix, which is stated as below.
In this section we prove the first statement \eqref{e:moment_converge} of Theorem \ref{semicircle}, and the following theorem, which states that the moments of the empirical eigenvalue distribution of $M^{(n)}$ concentrate on its expectation. 

\begin{theorem}
    \label{concentration}
   Consider the degree sequence $\{ d_1^{(n)}, d_2^{(n)}, \cdots, d_n^{(n)} \}$ for each $n\geq 1$. Suppose that $\min\{d_i^{(n)}\}_{i\in [n]} \gg \sqrt{D^{(n)}/n}$ as $n \rightarrow \infty$. Consider the normalized Laplacian matrix $M^{(n)}$ as defined in \eqref{e:defM}. Then, the moments of the empirical eigenvalue distribution of $M^{(n)}$ concentrates on its expectation: for any $\varepsilon > 0$, 
    $$ \lim_{n\rightarrow \infty} \p \biggl( \Big|\frac{1}{n}\Tr[({M^{(n)}})^k] - \frac{1}{n}\E[\Tr[({M^{(n)}})^k]] \Big| > \varepsilon \biggr) \to 0.
    $$
\end{theorem}
We use subsections \ref{subsec2.1}, \ref{subsec2.2}, \ref{subsec2.3}, and \ref{subsec2.4} to prove \eqref{e:moment_converge} of Theorem \ref{semicircle} and use subsection \ref{subsec2.5} to prove Theorem \ref{concentration}.

The second statement \eqref{e:weak_converge} of Theorem \ref{semicircle} follows from combining \eqref{e:moment_converge} and Theorem \ref{concentration} by a standard argument. We thus postpone its proof to the Appendix \ref{sec:appendix1}.

\subsection{Moments Expansion}\label{subsec2.1}
In this section we derive an expression for $\Tr[M^k]/n$. Let 
\begin{align}
    Q = \sqrt{D}\Delta^{-\frac{1}{2}}\left(A-\left[\frac{d_id_j}{D-1}\right]\right)\Delta^{-\frac{1}{2}} = n^\frac{1}{2}M,\quad Q_{ij} = \sqrt\frac{D}{d_id_j}\left(A_{ij}-\left[\frac{d_id_j}{D-1}\right]\right).
\end{align}
With the above notation, we have 
\begin{align}\begin{split}\label{e:Mk_momment}
    &\phantom{{}={}}\frac{1}{n}\E[\Tr(M^k)] =\frac{1}{n^{1+\frac{k}{2}}}\E\left[\Tr (Q^k)\right]= \sum_{i_1,i_2,\cdots ,i_k}\E\left[Q_{i_1i_2}Q_{i_2i_3}\cdots Q_{i_ki_1}\right]\\
     &= \frac{D^{\frac{k}{2}}}{n^{\frac{k}{2}+1}}
    \sum_{i_1,i_2,\cdots ,i_k}
    \frac{1}{d_{i_1}d_{i_2}\cdots d_{i_k}}\E\biggl[\left(A_{i_1i_2}-\left[\frac{d_{i_1}d_{i_2}}{D-1}\right]\right)\cdots\left(A_{i_ki_1}-\left[\frac{d_{i_k}d_{i_1}}{D-1}\right]\right)\biggr].
\end{split}\end{align}

We recall the indicator function $H_{ij}^{st}$ from \eqref{e:indicator_function}, which equals $1$ if the $s$-th half-edge of $v_i$ is connected to the $t$-th half-edge of $v_j$. Then, we have that $A_{ij} - \frac{d_id_j}{D-1}=\sum_{s=1}^{d_{i}}\sum_{t=1}^{d_{j}} (H_{ij}^{st} - \frac{1}{D-1})$.
And it follows that
\begin{align*}
    &\E\left[\left(A_{i_1i_2}-\left[\frac{d_{i_1}d_{i_2}}{D-1}\right]\right)\left(A_{i_2i_3}-\left[\frac{d_{i_2}d_{i_3}}{D-1}\right]\right)\cdots \left(A_{i_ki_1}-\left[\frac{d_{i_k}d_{i_1}}{D-1}\right]\right)\right]\\
    =&\sum_{s_1=1}^{d_{i_1}}\sum_{t_2=1}^{d_{i_2}}\,\sum_{s_2=1}^{d_{i_2}}\sum_{t_3=1}^{d_{i_3}}\cdots \sum_{s_k=1}^{d_{i_k}}\sum_{t_1=1}^{d_{i_1}} \E\biggl[\left(H_{i_1i_2}^{s_1t_2} - \frac{1}{D-1}\right)\left(H_{i_2i_3}^{s_2t_3} - \frac{1}{D-1}\right)\cdots \biggl(H_{i_ki_1}^{s_kt_1} - \frac{1}{D-1}\biggr)\biggr].
\end{align*}

In the rest of this section, we will need to evaluation
\begin{align}\label{e:product}
    \E\biggl[\left(H_{i_1i_2}^{s_1t_2} - \frac{1}{D-1}\right)\left(H_{i_2i_3}^{s_2t_3} - \frac{1}{D-1}\right)\cdots \biggl(H_{i_ki_1}^{s_kt_1} - \frac{1}{D-1}\biggr)\biggr]
\end{align}
which involves the most important definition as below.
\begin{definition} \label{def:kwalk}
    For any term in the form of \eqref{e:product}, we view the term as a closed path with vertices $v_{i_1}, v_{i_2}, v_{i_3}, \cdots, v_{i_{k}}, v_{i_1}$, where for each $1 \leq j \leq k$, the edge from $v_{i_j}$ to $v_{i_{j+1}}$ involves the $s_j$-th half-edge of vertex $v_{i_j}$ and the $t_{j+1}$-th half-edge of vertex $v_{i_{j+1}}$. We define such a closed path of length $k$ as a \textbf{$\textbf{k}$-walk}. 
\end{definition}
In general, the indicator functions $H_{i_1i_2}^{s_1t_2}, H_{i_2i_3}^{s_2t_3}, \ldots, H_{i_ki_1}^{s_kt_1}$ are weakly correlated, except when they share some half-edges. For example, if $t \neq t'$ or $j \neq j'$, at most one of $H_{ij}^{st}$ and $H_{ij'}^{st'}$ can be equal to $1$. To deal with such cases, we break the walk into two parts

% which can be viewed as a closed path with vertices $v_{i_1}, v_{i_2}, v_{i_3}, \cdots, v_{i_{k}}, v_{i_1}$, where for each $1 \leq j \leq k$, the edge from $v_{i_j}$ to $v_{i_{j+1}}$ involves the $s_j$-th half-edge of vertex $v_{i_j}$ and the $t_{j+1}$-th half-edge of vertex $v_{i_{j+1}}$. We define such a closed path as a $k$-walk. In general, the indicator functions $H_{i_1i_2}^{s_1t_2}, H_{i_2i_3}^{s_2t_3}, \ldots, H_{i_ki_1}^{s_kt_1}$ are weakly correlated, except when they share some half-edges. For example, if $t \neq t'$ or $j \neq j'$, at most one of $H_{ij}^{st}$ and $H_{ij'}^{st'}$ can be equal to $1$. To deal with such cases, we break the walk into two parts
% \begin{figure}
%     \centering
%     \includegraphics[width=0.5\linewidth]{Example1.png}
%     \caption{An example of ``good'' (in green) and ``bad'' (in red) edges. }
%     \label{fig:example1}
% \end{figure}

\begin{figure}
    \begin{minipage}[c]{0.49\textwidth}
        \includegraphics[width=\linewidth]{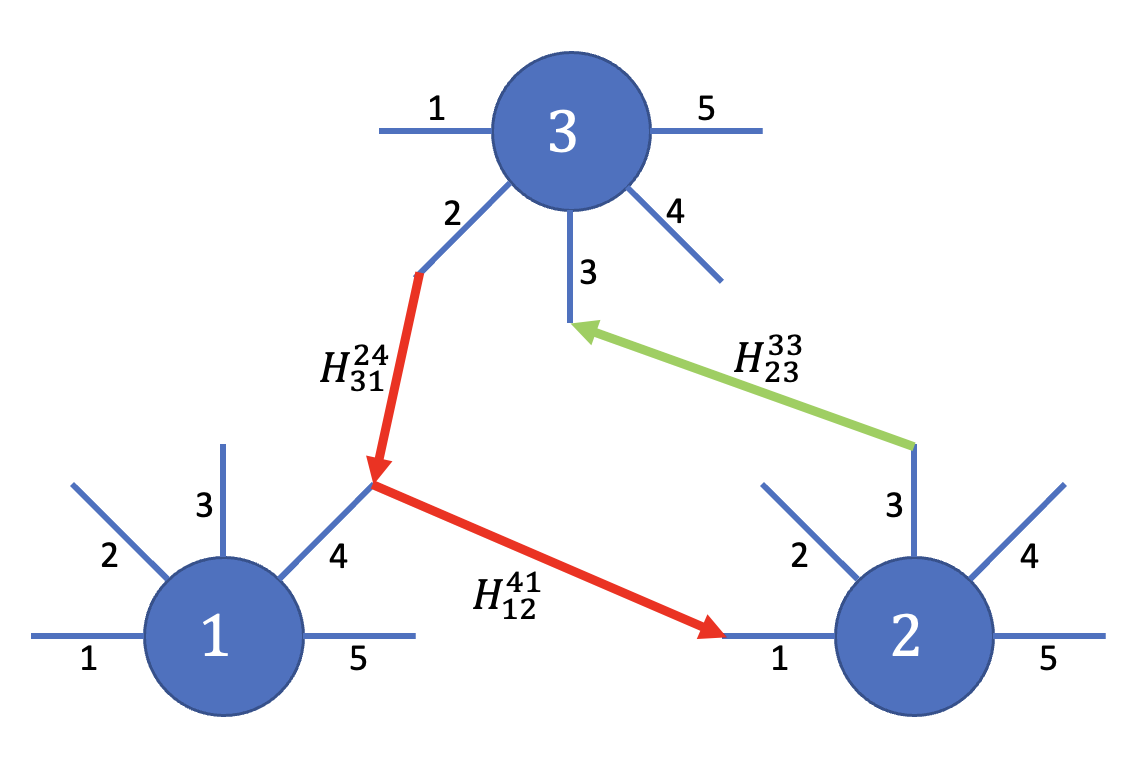}
        % \caption{Spectral distribution \\ for $d_1=100$, $d_2=500$, $n=10^5$.}
    \end{minipage}
    \hfill
    \begin{minipage}[c]{0.49\textwidth}
        \includegraphics[width=\linewidth]{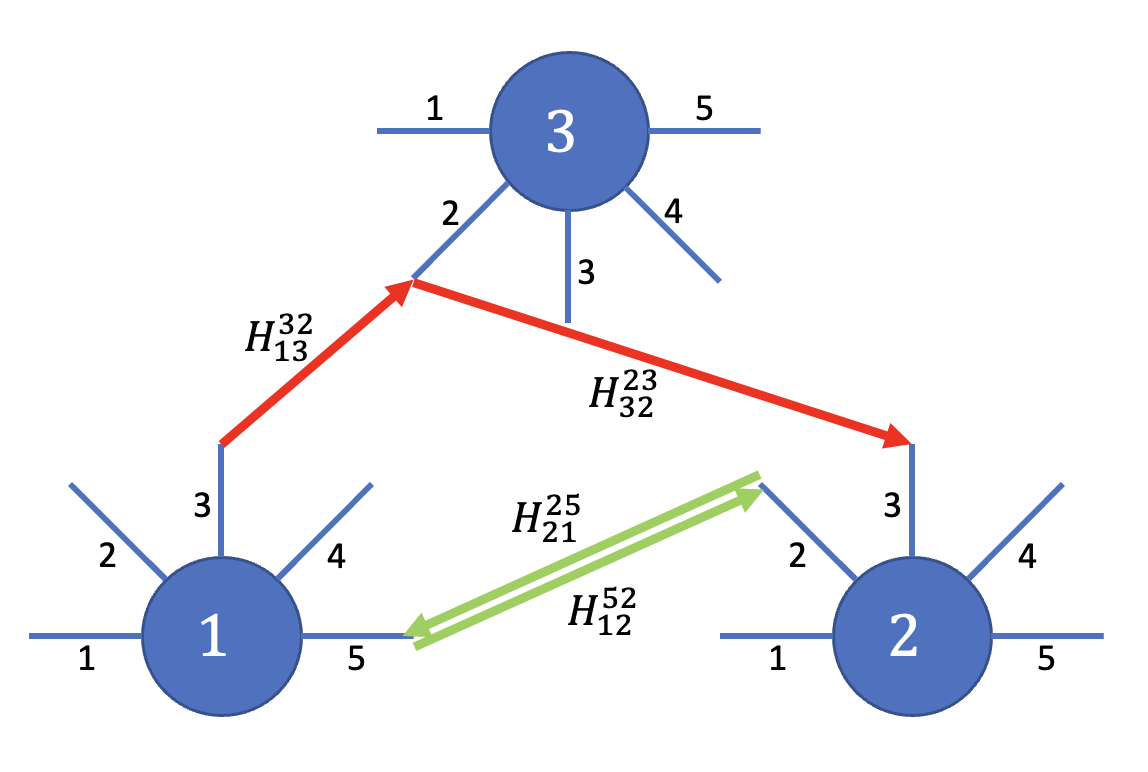}
        % \caption{Spectral distribution \\
        % for $d_1=10$, $d_2=200$, $n=10^5$.}
    \end{minipage}%
    \caption{Two examples of walks with ``good'' (in green) and ``bad'' (in red) edges.}
    \label{fig:k-walks}
\end{figure}

% \begin{figure}
%     \centering
%     \includegraphics[width=0.5\linewidth]{Example2.png}
%     \caption{Enter Caption}
%     \label{fig:enter-label}
% \end{figure}
% \begin{figure}
%     \centering
%     \includegraphics[width=0.5\linewidth]{Example3.png}
%     \caption{Enter Caption}
%     \label{fig:enter-label}
% \end{figure}
\begin{itemize}
    \item (``bad" part) For any $1\leq j\leq k$, $H_{i_ji_{j+1}}^{s_jt_{j+1}}$ is in ``bad" part if it shares exactly one half edge with another term. More specifically, there exists some other $1\leq j\neq j'\leq k$, such that
    \begin{align}
        |\{(i_j, s_j), (i_{j+1}, t_{j+1})\}\cap \{(i_{j'}, s_{j'}), (i_{j'+1}, t_{j'+1})\}|= 1.
    \end{align}
    In this case we call the $j$-th step a ``bad" edge, which consists of  the $s_j$-th half-edge of vertex $v_{i_j}$ and the $t_{j+1}$-th half-edge of vertex $v_{i_{j+1}}$.
    \item (``good" part) For any $1\leq j\leq k$, $H_{i_ji_{j+1}}^{s_jt_{j+1}}$ is in ``good" part if it does not share exactly one half edge with other terms. More specifically, for any other $1\leq j\neq j'\leq k$, we have
    \begin{align}
        |\{(i_j, s_j), (i_{j+1}, t_{j+1})\}\cap \{(i_{j'}, s_{j'}), (i_{j'+1}, t_{j'+1})\}|\neq 1.
    \end{align}
    We remark that we allow that $H_{i_ji_{j+1}}^{s_jt_{j+1}}$ share both half edges with some other term. In this case, then two indicator functions equal. 
\end{itemize}
Some examples of ``good'' and ``bad'' edges are illustrated in Figure \ref{fig:k-walks}. We remark that the number of terms in ``bad" parts is either zero, or is at least two, because there cannot be a single ``bad" edge. Moreover, the number of walks with a ``bad" part should be very small because it's extremely rare to pick two edges that share a half-edge.
% {\color{red} [I reformulate the bad and good parts, can you check it.]} {\color{blue} [I just checked it and added two examples.]} 

%Now let us characterize such a walk. The walk can be broken up into two parts: the ``good" part and the ``bad" part. The ``good" part denotes the portion of the walk where there are no two edges that share the same half-edge. Note that it is allowable for a ``good" part to contain two edges that are entirely the same. That is, two edges may use the same pair of half-edges coming from the same pair of vertices. We are ruling out the cases where two edges share exactly one half-edge. The ``bad" part denotes the complement portion of the walk, where each consecutive pair of edges shares the same half-edge from a vertex. For example, if we (luckily) pick $t_2$ and $s_2$ as the same half-edge from the vertex $i_2$, then we say that the edge ($s_1$, $t_2$) and the edge ($s_2$, $t_3$) are a pair of ``bad" edges. Note that ``bad" edges appear at least in the number of two because there cannot be a single ``bad" edge. Note that the number of walks with a ``bad" part should be very small because it's extremely rare to pick two edges that share a half-edge.

% We denote the edges in the ``bad" part as ``bad" edges. 

Let us rename the indicator random variables. For the  ``good" part, assume there are $r$ edges that appear multiple times in the product \eqref{e:product}, which we rename as $H_1, H_2, \dots, H_r$. Further, assume that each $H_i$ appears $\alpha_i \geq 2$ times.
For the remaining ``good" part, there are $p$ distinct edges, each appearing only once in the product \eqref{e:product}. We denote the corresponding indicators as $H_{r+1}, H_{r+2}, \dots, H_{r+p}$.
For the ``bad" part, assume there are $L$ distinct ``bad" edges, denoted as $H_{r+p+1}, H_{r+p+2}, \dots, H_{r+p+L}$, where each $H_i$ appears $\beta_i \geq 1$ times for $r+p+1 \leq i \leq r+p+L$. After renaming, we can rewrite \eqref{e:product} as follows:
% \begin{align*}
%     &\E\left[\left(H_{i_1i_2}^{s_1t_2} - \frac{1}{D-1}\right)\left(H_{i_2i_3}^{s_2t_3} - \frac{1}{D-1}\right)\cdots \left(H_{i_ki_1}^{s_kt_1} - \frac{1}{D-1}\right)\right] \\
%     &= \E\biggl[\left(H_1-\frac{1}{D-1}\right)^{\alpha_1}\left(H_2-\frac{1}{D-1}\right)^{\alpha_2}\cdots \left(H_r-\frac{1}{D-1}\right)^{\alpha_r}\\&\left(H_{r+1}-\frac{1}{D-1}\right)\cdots \left(H_{r+p}-\frac{1}{D-1}\right)
%     \\&\left(H_{r+p+1} - \frac{1}{D-1}\right)^{\beta_{r+p+1}}\cdots \left(H_{r+p+L_1} - \frac{1}{D-1}\right)^{\beta_{r+p+L_1}}
%     \\&\left(H_{r+p+L_1+1} - \frac{1}{D-1}\right)^{\beta_{r+p+L_1+1}}\cdots \left(H_{r+p+L_1+L_2} - \frac{1}{D-1}\right)^{\beta_{r+p+L_1+L_2}}
%     \\&\cdots 
%     \\&\left(H_{r+p+L_1+\cdots +L_{C-1}+1} - \frac{1}{D-1}\right)^{\beta_{r+p+L_1+\cdots +L_{C-1}+1}}\cdots \left(H_{r+p+L_1+\cdots +L_C} - \frac{1}{D-1}\right)^{\beta_{r+p+L_1+\cdots +L_C}}
%     \biggr]
% \end{align*}
\begin{align}\begin{split}\label{e:product2}
    &\E\left[\left(H_{i_1i_2}^{s_1t_2} - \frac{1}{D-1}\right)\left(H_{i_2i_3}^{s_2t_3} - \frac{1}{D-1}\right)\cdots \left(H_{i_ki_1}^{s_kt_1} - \frac{1}{D-1}\right)\right] \\
    = &\E\biggl[\left(H_1-\frac{1}{D-1}\right)^{\alpha_1}\left(H_2-\frac{1}{D-1}\right)^{\alpha_2}\cdots \left(H_r-\frac{1}{D-1}\right)^{\alpha_r}
    \\
    \times&\left(H_{r+1}-\frac{1}{D-1}\right)\cdots \left(H_{r+p}-\frac{1}{D-1}\right)
    \\
    \times &\left(H_{r+p+1} - \frac{1}{D-1}\right)^{\beta_{r+p+1}}\cdots \left(H_{r+p+L} - \frac{1}{D-1}\right)^{\beta_{r+p+L}}
    \biggr]
\end{split}\end{align}
where $\Sigma_{i=1}^{r}\alpha_i + p + \Sigma_{i=r+p+1}^{r+p+L} \beta_i = k$. Note that $\alpha_i \geq 2 \ \forall i$ and $\beta_i \geq 1 \ \forall i$.

We need to handle three  parts involved in the product \eqref{e:product2} differently. Recall that they are the ``good" part where each indicator random variable is repeated (denoted as \textbf{``good'' repeated edges}), the ``good" part where each indicator random variable has multiplicity $1$ (denoted as \textbf{``good' simple edges}), and the ``bad" part (denoted as \textbf{``bad'' edges}). For simplicity, when we are dealing with a part of the expectation, we usually rename the indicator random variables so that the index starts with $1$. For example, when we are later dealing with the ``good" simple edges, we use the notation of $H_1, \cdots , H_p$ instead of the cumbersome notation of $H_{r+1}, .., H_{r+p}$. 

\subsection{Bounds}\label{subsec2.2}
In this section, we respectively show upper bounds on the three parts in the product \eqref{e:product2} respectively.

First, we show the bounds for the ``good" simple part.
\begin{proposition} [Bounds on ``good" simple edges]
    \label{goodunique}
    Suppose there are $p$ distinct ``good" edges with multiplicities $1$. When conditioned on the first $r$ ``good" repeated edges and the last $L$ ``bad'' edges, the following bound is achieved
    $$\left|\E\left[\left(H_{r+1}-\frac{1}{D-1}\right)\cdots \left(H_{r+p}-\frac{1}{D-1}\right)\;\middle| \; H_1,H_2,\cdots ,H_r, H_{r+p+1}, \cdots , H_{r+p+L} \right] \right|\leq \frac{C_p}{(D-1)^\frac{3p}{2}} .$$
\end{proposition}

First note that after conditioning on $H_1,H_2,\cdots ,H_r$ and $H_{r+p+1}, \cdots , H_{r+p+L}$, the law of $H_{r+1}, H_{r+2},\cdots, H_{r+p}$ is still from a configuration model, which differs from the original configuration model only by a factor of $(1+o(1))$. This is because the new configuration model has degree sequence $d_i-O(k)=(1+o(1))d_i$.

Next we rename $\E[(H_{r+1}-\frac{1}{D-1})\cdots (H_{r+p}-\frac{1}{D-1})]$ as $\E[(X_1-\frac{1}{D-1})\cdots (X_p-\frac{1}{D-1})]$, and  bound its value using the techniques from \cite{Sa}.

Define the function $p_i:\{0,1\}^{i-1}\rightarrow \mathbb{R}$ as $p_i(x_1,x_2,\cdots ,x_{i-1})=\p(X_i=1|X_1=x_1,\cdots ,X_{i-1}=x_{i-1})$. Furthermore, for any $T = \{t_1, t_2,\cdots ,t_j\}
\subseteq [i-1]$, let $x_T = x_{t_1}x_{t_2}\cdots x_{t_j}$ and $p_i(T)$ be the function $p_i$ where $x_l = 1$ if $l\in T$ and $0$ otherwise. From \cite{OD}, we can write $p_i(x_1,\cdots ,x_{i-1})=\sum_{T\subseteq [i-1]}\chi_T(p_i)x_T$ for some unique $\chi_T(p_i)$. These are the Fourier coefficients of $p_i$. From \cite{OD}, $$\chi_T(p_i)=\sum_{S\subseteq T}(-1)^{|T|+|S|}p_i(S).$$

\begin{lemma} Suppose $p^3 \ll D$, $1\leq i\leq p$, and  $S\subseteq [i-1]$ is an arbitrary subset. Then, the following bound on the Fourier coefficient can be achieved
    \begin{align*}
      \chi_\emptyset(p_i)=\frac{1}{D-1},\quad   \left|\chi_S(p_i)\right| \leq \frac{(1+o(1))}{D-1} \left(\frac{2|S|}{D-1}\right)^{|S|} .
    \end{align*}
\end{lemma}

\begin{proof}
    The probability that $X_i=1$ is determined by the number of available half-edges. Each $X_j=1$ for $j<i$ makes two fewer available half-edges, so
    \begin{align}\label{e:piexp}
        p_i(X_i=1|X_1=x_1,\cdots ,X_{i-1}=x_{i-1}) &= \frac{(1+o(1))}{D-1-2\sum_{j=1}^{i-1}x_j} = \frac{(1+o(1))}{(D-1)\left(1-\frac{2\sum_{j=1}^{i-1}x_j}{D-1}\right)}. 
    \end{align}

% {\color{red}[Do we still need the following two paragraphs, since now we are restricted to good part?]}
% {\color{blue}[Sorry for the mistake! They are unnecessary since it's restricted to the good part. I have deleted them]}
    % It is also necessary to note here that some pairs of $X_k$ and $X_l$ could share half-edges. If $1 \leq k,l \leq i-1$, we assume in this proposition that at most one of $X_k$ and $X_l$ is $1$. We later show any case we consider satisfies this assumption.

    % In the case where $X_{i}$ shares half-edges with some $X_{j_1},X_{j_2},\cdots ,X_{j_y}$ where $1\leq j_1,j_2,\cdots ,j_y\leq i-1$, if any of $X_{j_1},\cdots ,X_{j_y}$ equal $1$, then $p_i=0$ since a half-edge cannot connect to two different half-edges. In this case, all the Fourier coefficients are identically $0$, so the inequality is trivially satisfied. If $X_{j_1}=\cdots =X_{j_y}=0$, there is an extra $y$ subtracted from the denominator of the above equation so that $p_i=\frac{1}{D-1-2\sum_{j=1}^{i-1}x_j-y}(1+o(1))$. The extra $y$ subtracted is due to the shared half-edge now having $y$ fewer options. This extra $y$ term can then be factored into the $(1+o(1))$ term.
    
    We can expand \eqref{e:piexp} as an infinite series,
    \begin{align*}
        \frac{1}{D-1}\frac{1}{1-\frac{2\sum_{j=1}^{i-1}x_j}{D-1}}&= \frac{1}{D-1}\sum_{l=0}^\infty \left(\frac{2\sum_{j=1}^{i-1}x_j}{D-1}\right)^l= \frac{1}{D-1}\sum_{l=0}^\infty \frac{2^l}{(D-1)^l}\left(\sum_{j=1}^{i-1}x_j\right)^l\\
        &= \frac{1}{D-1}\sum_{l=0}^\infty \frac{2^l}{(D-1)^l}\sum_{\alpha_1+\alpha_2+\cdots +\alpha_{i-1}=l}\frac{l!}{\prod_{j=1}^{i-1}\alpha_j!}\prod_{j=1}^{i-1}x_j^{\alpha_j}\\
        &= \frac{1}{D-1}\sum_{l=0}^\infty \frac{2^l}{(D-1)^l}\sum_{\alpha_1+\alpha_2+\cdots +\alpha_{i-1}=l}\frac{l!}{\prod_{j=1}^{i-1}\alpha_j!} \prod_{j=1}^{i-1}x_j^{\min\{1,\alpha_j\}}.
       \end{align*}
       From \cite{OD}, the Fourier coefficients are unique, so we have that the Fourier coefficient of $S=\{j_1, j_2,\cdots, j_m\}$ where $|S|=m$ is
        \begin{align}\label{e:chiSPi}
        \chi_S(p_i) &= \frac{1}{D-1} \sum_{l=m}^\infty \frac{2^l}{(D-1)^l} \sum_{\substack{\alpha_{j_1} + \alpha_{j_2} + \cdots  + \alpha_{j_m} = l\\ \alpha_{j_1}, \alpha_{j_2}, \cdots , \alpha_{j_m} \geq 1}} \frac{l!}{\prod_{q=1}^m \alpha_{j_q}!m!}(1+o(1))
    \end{align}
    where the sum
    \begin{align}\label{e:Ssum}
    \sum_{\substack{\alpha_{j_1} + \alpha_{j_2} + \cdots  + \alpha_{j_m} = l\\ \alpha_{j_1}, \alpha_{j_2}, \cdots , \alpha_{j_m} \geq 1}} \frac{l!}{\prod_{q=1}^m \alpha_{j_q}!}
    \end{align}
    is the number of ways to arrange $l$ marbles from $m$ different colors in a line, with at least 1 from each color. This is equivalent to the number of ways to put $l$ distinct marbles into $m$ boxes, with no boxes empty. This is $S(l,m)m!$, the Stirling numbers of the second kind. Therefore, $\eqref{e:Ssum}/m!
    = S(l, m)$. We can therefore provide an upper bound where the boxes can be empty so that $\eqref{e:Ssum}/m! \leq \frac{m^l}{m!}$.

    From here, we can bound \eqref{e:chiSPi} as
    \begin{align*}
        \chi_S(p_i) &\leq \frac{(1+o(1))}{D-1} \sum_{l=m}^\infty \frac{(2m)^l}{(D-1)^l}= \frac{(1+o(1))}{D-1}\frac{\left(\frac{2|S|}{D-1}\right)^{|S|}}{1-\frac{2|S|}{D-1}}= \frac{(1+o(1))}{D-1}\left(\frac{2|S|}{D-1}\right)^{|S|}
    \end{align*}
    which concludes the proof.
\end{proof}

\begin{lemma}[Proposition 3.2 from \cite{Sa}] Suppose there are $p$ ``good" simple edges. Then, their contribution can be evaluated using the following Fourier coefficients
\label{fourier}
    \begin{align*}
        &\phantom{{}={}}\E\left[\left(X_1-\frac{1}{D-1}\right)\left(X_1-\frac{1}{D-1}\right)\cdots\left(X_p-\frac{1}{D-1}\right)\right]\\
        =&\sum_{\substack{S_1,\cdots,S_p\\S_i\subseteq[i-1]}}\prod_{i=1}^p\left(1-\frac{1}{D-1} \cdot \mathds{1}_{i\in\cup_{j>i}S_j}\right)\chi_{S_i}\left(p_i-\frac{1}{D-1}\cdot \mathds{1}_{i \notin \cup_{j>i}S_j}\right).
    \end{align*}
\end{lemma}
The proof is omitted, since it is the same as in \cite{Sa}. 

\begin{proof}[Proof of \Cref{goodunique}]
Using Lemma \ref{fourier}, we can further bound the expectation as 
\begin{align} \label{e:chiprod}
    \E\left[\left(X_1-\frac{1}{D-1}\right)\cdots\left(X_p-\frac{1}{D-1}\right)\right]&\leq \sum_{\substack{S_1,\cdots,S_p\\S_i\subseteq[i-1]}}\prod_{i=1}^p\left|\chi_{S_i}\left(p_i-\frac{1}{D-1}\cdot \mathds{1}_{i \notin \cup_{j>i}S_j}\right)\right|.
\end{align}
From here, our proof strategy is to recreate the method in \cite{Sa}. For any set $S=(S_1,S_2,\cdots,S_p) \in \mathcal{D} = 2^{\emptyset}\times 2^{[1]} \times \cdots \times 2^{[p-1]}$, denote $\psi_S = \prod_{i=1}^p|\chi_{S_i}(p_i)|$. We iteratively define sets $\mathcal{C}_i$ and $\mathcal{D}_i$ and corresponding functions $f_i: \mathcal{D}_i\setminus \mathcal{C}_i \rightarrow \mathcal{D}_i$ such that $M\psi_S \leq \psi_{f_i(S)}$ for all $S \in \mathcal{D}_i\setminus \mathcal{C}_i$. Furthermore, if there are at most $N$ preimages for any $f_i(S)$, then we require $N/M = o(1)$. To see the motivation behind this, observe that
\begin{align*}
    \sum_{S \in \mathcal{D}_i} \psi_S &= \sum_{S\in \mathcal{D}_i\setminus \mathcal{C}_i} \psi_S + \sum_{S \in \mathcal{C}_i} \psi_S\\
    &\leq \sum_{S \in \mathcal{D}_i\setminus \mathcal{C}_i}\frac{1}{M}\psi_{f_i(S)} + \sum_{S \in \mathcal{C}_i}\psi_S\\
    &\leq \sum_{U \in f_i(\mathcal{D}_i\setminus \mathcal{C}_i)} \frac{N}{M} \psi_{U} + \sum_{S \in \mathcal{C}_i}\psi_S.
\end{align*}
Rearranging, we have
\begin{align*}
    \sum_{S \in \mathcal{C}_i} \psi_S &\geq (1-o(1))\sum_{S\in \mathcal{D}_i}\psi_S,
\end{align*}
so that the hypersets in $\mathcal{C}_i$ contribute most of the sum that makes up the Fourier coefficient.

We start with $\mathcal{D}_1 = \mathcal{D}$ and $\mathcal{C}_1=\{(S_1,\cdots,S_p) \in \mathcal{D}_1: S_i \cap S_j = \emptyset\}$. For any $U \in \mathcal{D}_1 \setminus \mathcal{C}_1$, there must be an element $a \in S_i \cap S_j$. We set $f_1:\mathcal{D}_1 \setminus \mathcal{C}_1 \rightarrow \mathcal{D}_1$ to be a function that takes the smallest $a\in S_i\cap S_j$ and removes $a$ from $S_i$ where $i < j$ without loss of generality. This multiplies the $i$-th bound in the product \eqref{e:chiprod} by $M \geq O(\frac{D-1}{p})(1+o(1))$. We also have that $N \leq p^2$ since there are at most $p$ choices for $a$ and $p$ choices for $S_i$. Then, as long as $p^3 \ll D$, $1-o(1)$ of the Fourier coefficient comes from the set $\mathcal{C}_1$. 

Next, we consider $\mathcal{D}_2 = \mathcal{C}_1$ and $\mathcal{C}_2=\{(S_1,\cdots,S_p) \in \mathcal{D}_2: \cup_{j=1}^p S_j \subseteq \{i: S_i=\emptyset\}\}$. For any $U \in \mathcal{D}_2 \setminus \mathcal{C}_2$, there must be $i \in \cup_{j=1}^p S_j$ such that $S_i \neq \emptyset $. We set $f_2: \mathcal{D}_2 \setminus \mathcal{C}_2 \rightarrow \mathcal{D}_1$ to be a function that considers the smallest-indexed such $S_i$ and removes $i$ from a later set $S_j$. This results in the same $M$ and $N$ bounds as before, so as long as $p^3 \ll D$, $1-o(1)$ of the Fourier coefficient comes from the set $\mathcal{C}_2$. 

We claim at least half of the sets in any $U \in \mathcal{C}_2$ must be empty. Suppose this is not the case. Then, fewer than half the sets are empty. All the other sets must take elements from these sets, so by the Pigeonhole principle, there must be at least two sets that share an element, which is a contradiction due to our definition of $\mathcal{C}_1$.

Next, we consider $\mathcal{D}_3 = \mathcal{C}_2$ and $\mathcal{C}_3 = \{(S_1,\cdots,S_p) \in \mathcal{D}_3: |S_i|\leq 1\}$. For any $U \in \mathcal{D}_3 \setminus \mathcal{C}_3$, there must be a set $S_i$ that contains at least two elements. We set $f_3: \mathcal{D}_3 \setminus \mathcal{C}_3 \rightarrow \mathcal{D}_3$ which takes the smallest-indexed set $S_i$ with at least two elements $a<b$, removes $a$ and $b$ from $S_i$, and adds $a$ to $S_b$. $S_i$'s bound changes by a factor of $O(\frac{(D-1)^2}{p^2})(1+o(1))$ while $S_b$'s bound changes by a factor of $O(\frac{p}{D-1})(1+o(1))$ so that $M \geq O(\frac{D-1}{p})(1+o(1))$. There are $p$ choices for $a$ and $p$ choices for $S_i$ so $N \leq p^2$. Then, if $p^3 \ll D$, $1-o(1)$ of the Fourier coefficient comes from the set $\mathcal{C}_3$. 

Next, we consider $\mathcal{D}_4 =  \mathcal{C}_3$ and $\mathcal{C}_4=\{(S_1,\cdots,S_p) \in \mathcal{D}_4: \sum_{i=1}^p|S_i| = \lfloor \frac{p}{2} \rfloor\}$. For any $U\in \mathcal{D}_4 \setminus \mathcal{C}_4$, we must have $\sum_{i=1}^p|S_i| < \lfloor \frac{p}{2} \rfloor$ by the claim that at least half of the sets are empty. Then, there must be a pair of $i < j$ such that $S_i=S_j=\emptyset$ and $i,j\notin S_k$ for any $k$ because the inequality is strict. We set $f_4: \mathcal{D}_4 \setminus \mathcal{C}_4 \rightarrow \mathcal{D}_4$ to be a function that considers the smallest such pair $i < j$ and adds $i$ to $S_j$. In this case, note that the bound on $S_i$ is $\chi_{S_i}(p_i - \frac{1}{D-1})=0$, so that $\mathcal{C}_4$ contributes $1-o(1)$ of the Fourier coefficient.

Then, $1-o(1)$ of $\psi_S$ comes from elements where half the sets are empty and half the sets have one element. For any $S \in \mathcal{C}_4$, we must have
\begin{align*}
    \psi_S &\leq \prod_{i: S_i = \emptyset}|\chi_\emptyset (p_i)| \prod_{j: S_j = \{k\} \text{ for some $k$}}|\chi_{\{k\}}(p_j)|^{\frac{p}{2}} \\
    &=\left(\frac{1}{D-1}\right)^{\frac{p}{2}}\left(\frac{1}{(D-1)^2}\right)^{\frac{p}{2}}\\
    &=\frac{1}{(D-1)^{\frac{3}{2}p}}.
\end{align*}
To bound $\E \left[\left(X_1-\frac{1}{D-1}\right)\cdots\left(X_p-\frac{1}{D-1}\right)\right]$, we must enumerate the number of such $S$. We pair all the elements up and then choose which element in each pair is the non-empty one, giving $2^{\frac{p}{2}}(p-1)!!$. Since $p$ is fixed, asymptotically, $\E \left[\left(X_1-\frac{1}{D-1}\right)\cdots\left(X_p-\frac{1}{D-1}\right)\right] \leq \frac{C_p}{(D-1)^{\frac{3}{2}p}}$.

Therefore, we have finished the proof for Proposition \ref{goodunique}.
\end{proof}
Next, we show the bounds for the ``good" repeated edges.
\begin{proposition} [Bounds on ``good" repeated edges] \label{goodrepeated}Suppose there are $r$ distinct ``good" repeated edges with multiplicities $\alpha_i\geq 2$. The following bound is achieved
    \begin{align*}
    \E\left[\left|H_1-\frac{1}{D-1}\right|^{\alpha_1}\left|H_2-\frac{1}{D-1}\right|^{\alpha_2}\cdots\left|H_r-\frac{1}{D-1}\right|^{\alpha_r}\right]\leq\frac{ (1+o(1))}{(D-1)^r}.
    \end{align*}
\end{proposition}
\begin{proof}
    % We show the case via induction on $r$. The base case for $r=1$ and $\alpha_1=2$ is easily verified. For $r+1$, we have
    % \begin{align*}
    %     &\E\left[\left|H_1-\frac{1}{D-1}\right|^{\alpha_1}\left|H_2-\frac{1}{D-1}\right|^{\alpha_2}\cdots\left|H_{r+1}-\frac{1}{D-1}\right|^{\alpha_{r+1}}\right]\\
    %     &= \E\left[\E\left[\left|H_{r+1}-\frac{1}{D-1}\right|^{\alpha_{r+1}}\;\middle|\; H_1,H_2,\cdots,H_r\right]|X|\right]\\
    %     &\leq \E\left[\E\left[\left|H_{r+1}-\frac{1}{D-1}\right|^{2}\;\middle|\; H_1,H_2,\cdots,H_r\right]|X|\right]\\
    %     &\leq \E\left[\E\left[\left(\frac{1}{D-1-2r}\right)^2-\frac{1}{(D-1)^2}\;\middle|\; H_1,H_2,\cdots,H_r\right]|X|\right]\\
    %     &\leq \frac{1}{(D-1)^{r+1}}
    % \end{align*}
    We show the case via induction on $r$. The base case for $r=1$ and $\alpha_1=2$ is easily verified. Now we assume the inductive hypothesis for $r = k$ that
    \begin{align*}
    \E\left[\left|H_1-\frac{1}{D-1}\right|^{\alpha_1}\left|H_2-\frac{1}{D-1}\right|^{\alpha_2}\cdots\left|H_k-\frac{1}{D-1}\right|^{\alpha_k}\right]\leq\frac{(1+o(1))}{(D-1)^k} 
    \end{align*}
    For $r = k+1$, we have
    \begin{align}\begin{split}
        &\phantom{{}={}} \E\left[\left|H_1-\frac{1}{D-1}\right|^{\alpha_1}\left|H_2-\frac{1}{D-1}\right|^{\alpha_2}\cdots\left|H_{k+1}-\frac{1}{D-1}\right|^{\alpha_{k+1}}\right]\\
        &= \E\left[\E\left[\left|H_{k+1}-\frac{1}{D-1}\right|^{\alpha_{k+1}}\;\middle|\; H_1,H_2,\cdots,H_k\right]|X|\right]\\
        &\leq \E\left[\E\left[\left|H_{k+1}-\frac{1}{D-1}\right|^{2}\;\middle|\; H_1,H_2,\cdots,H_k\right]|X|\right] \label{e:repeatedeq}
\end{split}    \end{align}    
    where 
    \begin{align}
        X:=\left|H_1-\frac{1}{D-1}\right|^{\alpha_1}\left|H_2-\frac{1}{D-1}\right|^{\alpha_2}\cdots\left|H_{k}-\frac{1}{D-1}\right|^{\alpha_{k}}.
    \end{align}
    The conditional expectation of $|H_{k+1}-1/(D+1)|^2$ is given by $(1+\oo(1))/(D-1)$, thus 
       \begin{align*} 
         &\phantom{{}={}} \E\left[\left|H_1-\frac{1}{D-1}\right|^{\alpha_1}\left|H_2-\frac{1}{D-1}\right|^{\alpha_2}\cdots\left|H_{k+1}-\frac{1}{D-1}\right|^{\alpha_{k+1}}\right]\\
        &\leq \E\left[ \frac{(1+o(1))}{D-1} |X|\right]\leq \frac{(1+o(1))}{D-1}  \cdot \frac{(1+o(1))}{(D-1)^k} \quad \text{ by the inductive hypothesis} \\
        &\leq \frac{(1+o(1))}{(D-1)^{k+1} }.
    \end{align*}
    \\
    Thus we conclude our proof.
\end{proof}
\begin{remark} \label{eqwhen2}
    The equality in Proposition \ref{goodrepeated} is achieved up to a factor of $1 + o(1)$ when $\alpha_i = 2$ for all $i\in [r]$. To see why, if all multiplicities are $2$, then the equality is achieved on the line of \ref{e:repeatedeq} and this is the only place where the inequality is used.
\end{remark}
Then, we show the bounds for the ``bad" part. 
\begin{proposition} [Bounds on ``bad" edges]
    \label{bad}
    Suppose there are $L$ distinct ``bad" edges with multiplicities $\beta_i\geq 1$. When conditioned on the first $r$ ``good" repeated edges, the following bound is achieved
    \begin{align*}
        &\E\left[\left|H_{r+p+1} - \frac{1}{D-1}\right|^{\beta_{r+p+1}}\cdots  \left|H_{r+p+L} - \frac{1}{D-1}\right|^{\beta_{r+p+L}}\middle|H_1, H_2, \cdots, H_{r} \right]\leq \frac{C_L}{(D-1)^{L}}.
    \end{align*}
\end{proposition}
\begin{proof}
    The same as in the proof of \Cref{goodunique}, we notice that conditioning on $H_1,H_2,\cdots,H_{r}$ only adds a factor of $(1+o(1))$. We only need to consider $\E[|H_{r+p+1}-\frac{1}{D-1}|^{\beta_{r+p+1}}\cdots|H_{r+p+L}-\frac{1}{D-1}|^{\beta_{r+p+L}}]$ which we rename as $\E[|X_1-\frac{1}{D-1}|^{\beta_1}\cdots|X_{L}-\frac{1}{D-1}|^{\beta_{L}}]$.
    
    We can encode $(X_1, X_2, \cdots, X_{L})$ by a subset $S\subseteq [L]$, where $i\in S$ if $X_i=1$, otherwise $X_i=0$. We say $S$ is realizable if it is possible to construct a corresponding graph with all edges in $S$. For example, $S$ is not realizable if there are two edges in $S$ which connect to the same half-edge of a vertex. Then, by using the definition of expectation, we have
    \begin{align*}
        &\E\left[\left|X_1-\frac{1}{D-1}\right|^{\beta_1}\cdots\left|X_{L}-\frac{1}{D-1}\right|^{\beta_{L}}\right]\\
        &= \left(1+O\left(\frac{1}{D}\right)\right) \sum_{S\subseteq [L]} \mathbbm{1}_{\text{S is realizable}} \left|\frac{1}{D-1}\right|^{\left|S\right|} \left|1 - \frac{1}{D-1}\right|^{L-\left|S\right|} \left|1 - \frac{1}{D-1}\right|^{\sum_{i\in S} \beta_i} \left|0 - \frac{1}{D-1}\right|^{\sum_{i\notin S} \beta_i}\\
        &\leq \left(1+O\left(\frac{1}{D}\right)\right) \sum_{S\subseteq [L]} \mathbbm{1}_{\text{S is realizable}} \left|\frac{1}{D-1}\right|^{|S|+\sum_{i\notin S} \beta_i} .
    \end{align*}
    The expression above is a sum of terms of different orders. Note that for any subset $S \subseteq [L]$, we must have $|S| + \sum_{i\notin S} \beta_i = \sum_{i\in S} 1 + \sum_{i \notin S} \beta_i \geq \sum_{i\in S} 1 + \sum_{i \notin S} 1 = L$. Now we know that each term in the above sum is bounded by $(D-1)^{-L}$. Since there are at most $2^L$ possible $S$'s, we know
    % Since $D$ goes to infinity, we only need to find out the dominant term. It is notable that the smallest order is when $S = [L]$. Though $S = [L]$ might not be realizable because the ``bad" edges correlate with each other, we observe that other choices of $S \subset [L]$ will only make the order larger. Thus, the term with $S = [L]$ can be a upper bound for the expression above. 
    \begin{align*}
        &\E\left[\left|X_1-\frac{1}{D-1}\right|^{\beta_1}\cdots\left|X_{L}-\frac{1}{D-1}\right|^{\beta_{L}}\right] \leq \frac{2^L}{(D-1)^{L}}\left(1+O\left(\frac{1}{D}\right) \right) 
        \leq \frac{C_L}{(D-1)^{L}}
    \end{align*}
    where $C_L$ is a constant only depending on $L$. Putting back the conditions, we get
    \begin{align*}
        &\E\left[\left|X_1-\frac{1}{D-1}\right|^{\beta_1}\cdots\left|X_{L}-\frac{1}{D-1}\right|^{\beta_{L}} \middle|H_1, H_2, \cdots, H_{r} \right]
        \leq \frac{C_L}{(D-1)^{L}} (1+o(1)).
    \end{align*}
    Thus we conclude our proof.
\end{proof}
Finally, we combine all the bounds together and prove the bound for a $k$-walk, as defined in Definition \ref{def:kwalk}, with parameters of $r, p, L$.
% {\color{red} Introduce the definition of $k$-walk, maybe at the beginning of this section.} {\color{blue} I introduced $k$-walk when we wrote the expression of the long product H.}
\begin{proposition} [Bounds on an arbitrary $k$-walk]
\label{p:k-walk_bound}
Suppose an arbitrary $k$-walk consists of $r$ ``good" repeated edges, $p$ ``good" simple edges, and $L$ ``bad" edges. Then the following bound is achieved
    \begin{align}\begin{split}\label{e:Hproduct}
        &\E\biggl[\left(H_1-\frac{1}{D-1}\right)^{\alpha_1}\left(H_2-\frac{1}{D-1}\right)^{\alpha_2}\cdots\left(H_r-\frac{1}{D-1}\right)^{\alpha_r}
        \\&\left(H_{r+1}-\frac{1}{D-1}\right)\cdots\left(H_{r+p}-\frac{1}{D-1}\right)
        \\&\left(H_{r+p+1} - \frac{1}{D-1}\right)^{\beta_{r+p+1}}\cdots\left(H_{r+p+L} - \frac{1}{D-1}\right)^{\beta_{r+p+L}}
        % \\&\left(H_{r+p+L_1+1} - \frac{1}{D-1}\right)^{\beta_{r+p+L_1+1}}\cdots\left(H_{r+p+L_1+L_2} - \frac{1}{D-1}\right)^{\beta_{r+p+L_1+L_2}}
        % \\&\cdots
        % \\&\left(H_{r+p+L_1+\cdots+L_{C-1}+1} - \frac{1}{D-1}\right)^{\beta_{r+p+L_1+\cdots+L_{C-1}+1}}\cdots\left(H_{r+p+L_1+\cdots+L_C} - \frac{1}{D-1}\right)^{\beta_{r+p+L_1+\cdots+L_C}}
        \biggr]
        \\&\leq \frac{(1+o(1))}{(D-1)^r} \frac{C_p}{(D-1)^{\frac{3p}{2}}} \frac{C_L}{(D-1)^L}.
    \end{split}\end{align}
    % where $L = L_1 + L_2 + \cdots + L_C$ is the number of distinct edges in the ``bad" part.
\end{proposition}

\begin{proof}
    Let $X$ denote the product of the first $r$ terms and the last $L$ terms on the lefthand side of \eqref{e:Hproduct}. By using conditional expectation, we have
    \begin{align*}
        LHS
        &\leq \E \biggr[  \E\biggr[ \left(H_{r+1} - \frac{1}{D-1}\right)\cdots\left(H_{r+p} - \frac{1}{D-1}\right) \Big| X \biggr] \lvert X \rvert \biggr]\\
        &\leq  \E\biggl[\frac{C_p}{(D-1)^{\frac{3p}{2}}}  \lvert X \rvert \biggl] \text{ by Proposition \ref{goodunique}}\\
        &= \frac{C_p}{(D-1)^{\frac{3p}{2}}} \E[\lvert X \rvert].
    \end{align*}
    Then we observe the expectation of the remaining terms. Let $Y$ denote the product of the first $r$ terms on the lefthand side of \eqref{e:Hproduct}. By using conditional expectation again, we have
    \begin{align*}
        \E [ \lvert X \rvert ]
        %&\leq \E \biggr[
        %\left| \left(H_1-\frac{1}{D-1}\right)^{\alpha_1}\cdots\left(H_r-\frac{1}{D-1}\right)^{\alpha_r} \right| \left|H_{r+p+1} - \frac{1}{D-1}\right|^{\beta_{r+p+1}}\cdots\left|H_{r+p+L} - \frac{1}{D-1}\right|^{\beta_{r+p+L}}
        %\biggr]\\
        % &= \E \biggr[ \E \biggr[ \left|H_{r+p+1} - \frac{1}{D-1}\right|^{\beta_{r+p+1}}\cdots\left|H_{r+p+L} - \frac{1}{D-1}\right|^{\beta_{r+p+L}}
        % \Big| Y \biggr] Y
        % \biggr] \\
        &\leq \E \biggr[ \E \biggr[ \left|H_{r+p+1} - \frac{1}{D-1}\right|^{\beta_{r+p+1}}\cdots\left|H_{r+p+L} - \frac{1}{D-1}\right|^{\beta_{r+p+L}}
        \Big|  Y  \biggr] \lvert Y \rvert
        \biggr] \\
        &\leq \E \biggr[
        \frac{C_L}{(D-1)^L}\lvert Y \rvert
        \biggr] \text{ by Proposition \ref{bad}}\\
        &= \frac{C_L}{(D-1)^L}\E [\lvert Y \rvert].
    \end{align*}
    Next, we apply Proposition \ref{goodrepeated} and get
    \begin{align*}
        \E[\lvert Y \rvert] &\leq \E\biggr[\left|H_1-\frac{1}{D-1}\right|^{\alpha_1}\cdots\left|H_r-\frac{1}{D-1}\right|^{\alpha_r} \biggr]
        % &\leq \E\biggr[\left|H_1-\frac{1}{D-1}\right|^{\alpha_1}\cdots\left|H_r-\frac{1}{D-1}\right|^{\alpha_r} \biggr]\\
        \leq \frac{1+o(1)}{(D-1)^r}.
    \end{align*}
    Finally, we combine all estimates above together and get
    \begin{align*}
        LHS \leq \frac{C_p}{(D-1)^{\frac{3p}{2}}} \frac{C_L}{(D-1)^L} \frac{(1+o(1)) }{(D-1)^r} 
    \end{align*}
    which concludes the proof.
\end{proof}

From here, we can put it all together. Our strategy now involves representing the entire sum we have bounded thus far, and then once again finding functions to reduce the size of the set until we are only counting trees with $k/2$ edges.

Thanks to \Cref{p:k-walk_bound}, we can upper bound \eqref{e:Mk_momment} as
\begin{align*}
    \frac{1}{n}|\E[\Tr M^k]| 
     &\leq  \frac{1}{n} \sqrt\frac{D}{n}^k\sum_{\substack{i_1,\cdots,i_k\\s_1,t_2,\cdots,s_k,t_1}} \frac{1}{d_{i_1}\cdots d_{i_k}} \frac{(1+o(1)) C_p C_L}{(D-1)^r(D-1)^\frac{3p}{2}(D-1)^L}\\
      &= \frac{1}{n} \sqrt\frac{D}{n}^k\sum_{\substack{i_1,\cdots,i_k\\s_1,t_2,\cdots,s_k,t_1}} \frac{1}{d_{i_1}\cdots d_{i_k}} \frac{(1+o(1)) C_p C_L}{D^rD^\frac{3p}{2}D^L}\\
     &= \frac{D^{\frac{k}{2}}}{n^{\frac{k}{2}+1}} \sum_{G\in \mathcal{G}_k} \sum_{i_1, \cdots, i_k} \sum_{s_1, t_2, \cdots, s_k, t_1} \frac{1}{d_{i_1}\cdots d_{i_k}} \frac{(1+o(1)) C_p C_L}{D^{r+\frac{3p}{2}+L}}
\end{align*}
where we take an approximation by replacing $D - 1$ with $D$. The set $\mathcal{G}_k$ denotes the set of all possible shapes of a $k$-walk up to permutations of vertices and half-edges. However, we distinguish ``good'' edges and ``bad'' edges in the structure. For example, a circle consisting of $k$ ``good'' edges is counted as a different shape than a circle consisting of $k$ ``bad'' edges.
% In the expression above, the first sum over $G$ determines the exact shape of the walk, without specifying the vertices or the half-edges used. For example, one such walk may look like $(a, b), (b, c), (c,a), (a,c), (c,a)$, where we don't know what $a, b, c$ stands for. It also determines whether there are ``bad" edges as part of the shape, such as $(c,a)$ and $(a, c)$ sharing a half-edge from $a$. 
% {\color{red}Introduce formally, what is the set of $G$ is summing over...}

In the expression above, the second sum over $i_1, \cdots, i_k$ determines what vertices are used to instantiate the walk. One such example may be $i_1 = 1, i_2 = n/2, i_3= n,\cdots$ out of the $n$ vertices that we consider. The third sum over $s_1, t_2,\cdots, s_k, t_1$ refers to the previous step that we break up the terms $A_{i_1i_2} - \frac{d_id_j}{D-1}$ into $\sum_{s_1=1}^{d_{i_1}}\sum_{t_2=1}^{d_{i_2}} (H_{i_1i_2}^{s_1t_2} - \frac{1}{D-1})$.

Now we focus on the following part
\begin{align*}
    \sum_{i_1, \cdots, i_k} \sum_{s_1, t_2, \cdots , s_k, t_1} \frac{1}{d_{i_1}\cdots d_{i_k}}.
\end{align*}

% \begin{align*}
%     \frac{1}{n}\E[\Tr M^k] &= \frac{1}{n} \sqrt\frac{D}{n}^k\sum_{\substack{i_1,\cdots ,i_k\\s_1,t_2,\cdots ,s_k,t_1}} \frac{1}{d_{i_1}\cdots d_{i_k}} \frac{1}{D^rD^\frac{3p}{2}}(1+o(1))
% \end{align*}

% Taking a closer look at $$\sum_{\substack{i_1,\cdots ,i_k\\s_1,t_2,\cdots ,s_k,t_1}} \frac{1}{d_{i_1}\cdots d_{i_k}}$$

We modify our notation a bit to better represent the results. Note that $i_1,\cdots ,i_k$ encode a $k$-walk and $s_1,t_2,\cdots ,s_k,t_1$ encode connected half-edges. We use $V$ to denote its vertex set.

Now let $E_k$ be all the edges traversed in the ``good" part of the $k$-walk and $V_k$ all the vertices involved the ``good" part. Let $E_k'$ be all the edges traversed in the ``bad" part of the $k$-walk and $V_k'$ all the vertices involved the ``bad" part. Note that $V_k$ and $V_k'$ may overlap and that $V = V_k \cup V_k'$.

For each ``good" edge $e\in E_k$ (including both repeated edges and simple edges), we denote its multiplicity as $\alpha_e$. For each ``bad" edge $e'\in E_k'$, we denote its multiplicity as $\beta_{e'}$. We also have that $\sum_e \alpha_e + \sum_{e'} \beta_{e'} = k$. For each vertex $i\in V_k$, we use $z_i$ to denote the number of distinct half-edges used for vertex $i$ that correspond to the ``good" edges. For each vertex $i\in V_k'$, we use $z_i'$ to denote the number of distinct half-edges used for vertex $i$ that correspond to the ``bad" edges. 

% Although not written above, each connected pair of half-edges $e$ is crossed $\alpha_e$ times. We have that $\sum_e \alpha_e + p = k$, with $q$ connected pairs of half-edges.

If we sum over all such graphs, we can write, for every ``good" edge $e = (i_e,j_e) \in E_k$ and every ``bad" edge $e' = (i_{e'}, j_{e'}) \in E_k'$,
\begin{align*}
    \sum_{i_1, \cdots , i_k} \sum_{s_1, t_2, \cdots , s_k, t_1} \frac{1}{d_{i_1}\cdots d_{i_k}} &= \sum_{i_1,\cdots ,i_k} \prod_{e \in E_k} \frac{1}{\sqrt{d_{i_e}d_{j_e}}^{\alpha_e}} \prod_{e \in E_k} d_{i_e} d_{j_e}
    \prod_{e' \in E_k'} \frac{1}{\sqrt{d_{i_{e'}}d_{j_{e'}}}^{\beta_{e'}}} \prod_{i \in V_k'} d_i^{z_i'}.
\end{align*}

Here is an important difference between the ``good" part and ``bad" part of the walk. For a ``good" edge $e\in E_k$, we need to divide by $\sqrt{d_{i_{e}}d_{j_{e}}}^{\alpha_{e}}$ to account for the multiplicity of $\alpha_e$. We further need to multiply by $d_{i_e} d_{j_e}$ to select two new half-edges, one from $i_e$ and one from $j_e$. It is not hard to verify that, for a vertex $i\in V_k$, $z_i = \sum_{i\in e, e\in E_k} 1$, since ``good" edges don't share half-edges. For a ``bad" edge $e' \in E_k'$, we still need to divide by $\sqrt{d_{i_{e'}}d_{j_{e'}}}^{\beta_{e'}}$ to account for the multiplicity of $\beta_{e'}$. However, when selecting half-edges, the definition of ``bad" part indicates that we may not need to select new half-edges for $i_{e'}$ and $j_{e'}$. Thus, we use $z_i'$ to denote the number of distinct half-edges used for vertex $i$ in the ``bad" part and loop over all vertices $i\in V_k'$ to account for selecting half-edges. 
% {\color{red}[it will be helpful to add a figure for these $E_k, E_k', z_i, z_i'$].}
\begin{figure}
    \centering
    \includegraphics[width=0.6\linewidth]{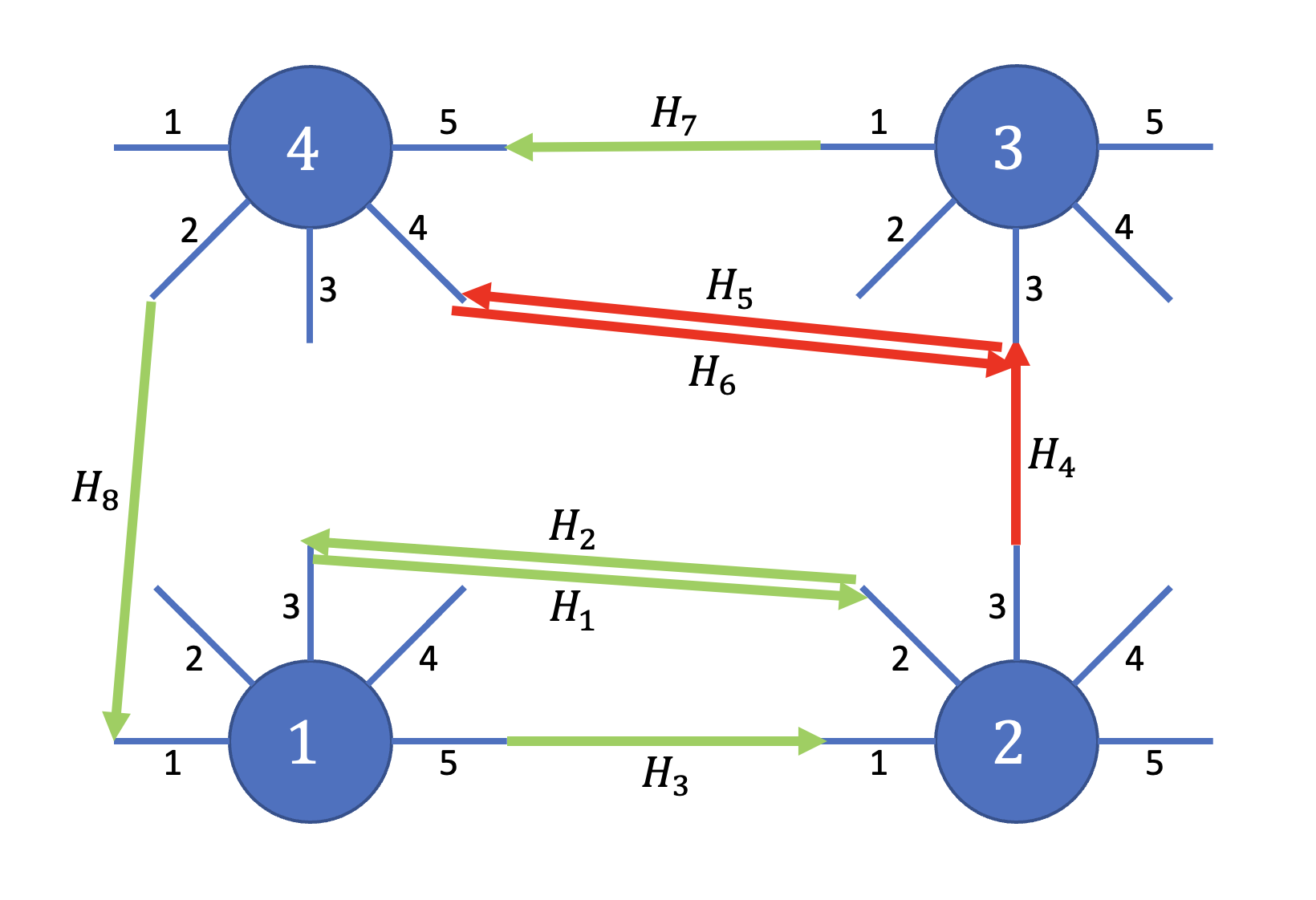}
    \caption{An example of a walk with ``good'' (in green) and ``bad'' (in red) edges.}
    \label{fig:exampleE_k}
\end{figure}

One example is illustrated in Figure \ref{fig:exampleE_k}, where our notations are defined in the following way,
\begin{align*}
    k = 8, \quad V_k = \{ 1, 2, 3, 4\}, &\quad V_k'= \{2,  3, 4\} \\
    E_k = \{ H_1, H_2, H_3, H_7, H_8\}, &\quad E_k'= \{ H_4, H_5, H_6\} \\
    z_1 = 3,\quad z_1' = 0,\quad z_2 = 2,\quad z_2' = 1,&\quad z_3 = 1,\quad z_3' = 1,\quad z_4 = 2,\quad z_4' = 1.
\end{align*}

Counting over vertices instead of edges, we have
\begin{align*}
    \sum_{i_1, \cdots , i_k} \sum_{s_1, t_2, \cdots , s_k, t_1} \frac{1}{d_{i_1}\cdots d_{i_k}} 
    % &= \sum_{i_1,\cdots ,i_k} \prod_{e \in E_k} d_{i_e}^{1 - \frac{\alpha_e}{2}} d_{j_e}^{1 - \frac{\alpha_e}{2}} \prod_{e' \in E_k} d_{i_{e'}}^{1 - \frac{\beta_{e'}}{2}} d_{j_{e'}}^{1 - \frac{\beta_{e'}}{2}}\\
    &= \sum_{i_1,\cdots ,i_k} \prod_{i \in V} d_i^{\sum_{i \in e, e\in E_k} (1 - \frac{\alpha_e}{2}) + z_i' - \sum_{i \in e', e'\in E_k'}  \frac{\beta_{e'}}{2}}\\
    &= \sum_{i_1,\cdots ,i_k} \prod_{i \in V} d_i^{ z_i - \sum_{i \in e, e\in E_k} \frac{\alpha_e}{2} + z_i' - \sum_{i \in e', e'\in E_k'}  \frac{\beta_{e'}}{2}}\\
    &= \sum_{i_1,\cdots ,i_k} \prod_{i \in V} d_i^{w_i}= \prod_{i \in V} \sum_{l = 1}^n d_l^{w_i}
\end{align*}
where we set
\begin{align}\label{e:defwi}
    w_i = z_i - \sum_{i \in e, e\in E_k} \frac{\alpha_e}{2} + z_i' - \sum_{i \in e', e'\in E_k'}  \frac{\beta_{e'}}{2}.
\end{align}
The last equation is because we want to specify all $i_1, \cdots , i_k$ with vertices from $1$ to $n$. By the Binomial Theorem, it's easy to establish the last equation. 

Summarizing what we have here, we establish that
\begin{align*}
    |\frac{1}{n} \E[\Tr M^k]|
    % &\leq  \frac{D^{\frac{k}{2}}}{n^{\frac{k}{2}+1}} \sum_{G\in \mathcal{G}_k} \sum_{i_1, \cdots, i_k} \sum_{s_1, t_2, \cdots, s_k, t_1} \frac{1}{d_{i_1}\cdots d_{i_k}} \frac{(1+o(1)) C_p C_L}{D^{r+\frac{3p}{2}+L}}
    % \\
    &\leq \sum_{G\in \mathcal{G}_k}  \frac{D^{\frac{k}{2}}}{n^{\frac{k}{2}+1}}
    \prod_{i \in V} \sum_{l = 1}^n d_l^{w_i}
    \frac{(1+o(1)) C_p C_L}{D^{r+\frac{3p}{2}+L}}\\
    &= \sum_{G\in \mathcal{G}_k}  \frac{(1+o(1)) C_p C_L}{n}\frac{D^{\frac{k}{2}}}{n^{\frac{k}{2}}}
    \prod_{i \in V} \sum_{l = 1}^n d_l^{w_i}
    \frac{1}{D^{r+\frac{3p}{2}+L}}\\
    &= \sum_{G\in \mathcal{G}_k}  \frac{(1+o(1)) C_p C_L}{n} t_G
\end{align*}
where we set 
\begin{align}\label{e:tG}
    t_G = \frac{D^{\frac{k}{2}}}{n^{\frac{k}{2}}} \prod_{i\in V}\sum_{l=1}^n d_l^{w_i}\frac{1}{D^{r+\frac{3p}{2}+L}}
\end{align}

% Next, we constrain what $G$ can look like to only Dyck paths by showing that any non-Dyck $G$ contributes $0$ asymptotically as $n\rightarrow \infty$. Note that $G$ currently consists of all graphs with all even degrees.

Before constraining the set of graphs, we first establish some important observations about the powers $w_i$, summarized in the following lemma. 
\begin{lemma} [Properties of $w_i$]
    \label{lemmawi}
    If $w_i = 1$, then $\sum_l d_l^{w_i} = D$. If $w_i =0$, then $\sum_l d_l^{w_i } = n$. If $w_i \leq -1$, then $\sum_l d_l^{w_i}  = o(n(n/D)^{|w_i|/2})$.\\
    Additionally, if $w_i$ increases by $1$, then $\sum_l d_l^{w_i}$ increases by a factor of $\omega (\sqrt{D/n})$. 
\end{lemma}
\begin{proof}
    The cases where $w_i = 1$ and where $w_i = 0$ are easy to verify. Now let us discuss the case $w_i \leq -1$. By the assumption of $d_i = \omega (\sqrt{D/n})$ for all $i$, we have that 
    $$\sum_l d_l^{w_i} = \sum_l d_l^{-|w_i|} = \sum_l o\left( \sqrt{\frac{n}{D}}^{|w_i|} \right) = o\left(n \sqrt{\frac{n}{D}}^{|w_i|}\right).$$
    Additionally, let us discuss the effect of increasing $w_i$ by $1$. Again by the assumption of $d_i = \omega (\sqrt{D/n})$ for all $i$, we know that each term of the sum increases by a factor of at least $\omega (\sqrt{D/n})$. Thus the sum will increase by a factor of at least $\omega (\sqrt{D/n})$.
\end{proof}

\begin{remark}
    The assumption ``$d_i \gg \sqrt{D/n}$ for all $i$" is crucial for the lemma to hold. If such assumption is substituted by the weaker condition ``$d_i \geq \sqrt{D/n}$ for all $i$", then we have a slightly weaker conclusion: if $w_i$ increases by $1$, $\sum_l d_l^{w_i}$ will increase by a factor of $\Omega (\sqrt{D/n})$. 
\end{remark}

% If $w_i  \leq -1$, then $\sum_l d_l^{w_i} = o( \frac{n}{\sqrt{\frac{D}{n}}^{|w_i| }} ) = o(n\sqrt{\frac{n}{D}}^{|w_i|})$, if $w_i =0$, then $\sum_l d_l^{w_i } = n$, and if $w_i = 1$, $\sum_l d_l^{w_i} = D$.

\subsection{Transformation}\label{subsec2.3}

We first consider any graph $G$ which contains at least one edge with odd multiplicity. Our plan here is to transform $G$ so that $t_{\widetilde G}$ of the final graph $\widetilde G$ will be convenient to deal with. 

In the analysis below, we will frequently but implicitly refer to the bound of $t_G$, whose expression is in \eqref{e:tG}. There are three parts: 
\begin{align}\label{e:tGcopy}
\frac{D^{\frac{k}{2}}}{n^{\frac{k}{2}}}, \quad \prod_{i\in V}\sum_{l=1}^n d_l^{w_i},  \quad \frac{1}{D^{r+\frac{3p}{2}+L}}.\end{align}
Whenever we make a transformation, we will comment the changes on these three parts respectively and summarize the overall change in the bound. Our goal is to make sure that each step of transformation will only increase the bound, so that as long as the bound of the final transformed graph goes to $0$ asymptotically, all such graphs have contribution of $0$ asymptotically.
\subsubsection{Transformation of the ``good" part }

Consider any $\alpha_e \geq 4$ for $e=(a,b)$. In these cases, we can see that we upper bound the expression's value by substituting $\alpha_e$ with $\alpha_e -2$. This increases $w_a$ by $1$ and $w_b$ in \eqref{e:defwi} by $1$, so the product increases by a factor of $\omega (D/n)$. Since $k$ decreases by $2$, the expression of $t_G$ in \eqref{e:tG} decreases by a factor of $n/D$. Thus, the overall bound increases by a factor of $\omega(1)$.

% Therefore, we only need to consider graphs where $\alpha_e = 1,2,3$. If $\alpha_e = 2$, the expression does not change if we substitute $\alpha_e \rightarrow \alpha_e - 2$. The reason is below. First of all, it is noteworthy that a small graph consisting only of $(a,b)$ with multiplicity $2$ has a contribution of $1$. The computation is quite simple,
% $$t_{(a,b)} = \frac{D^{\frac{2}{2}}}{n^{\frac{2}{2} + 1}} n^2 \frac{1}{D} = 1$$
% Thus, we can decompose the original graph $G$ into $(a,b)$ and $G'$. Note that $G'$ may be a connected graph or a disconnected graph, and either case works. In the latter case, we denote $G' = G_a \cup G_b$. By the fact that $\sum_{l=1}^n d_l^{w_i} \leq \sum_{l=1}^n d_l^{w_i'}\sum_{l=1}^n d_l^{w_i''}$ where $w_i'+w_i''=w_i$ for all $i$, we know that $t_{G} \leq t_{(a,b)} t_{G'}$. Thus, we can freely drop the edges with multiplicities $2$. In the disconnected case, we have $t_G \leq t_{(a,b)} t_{G_a} t_{G_b}$. We can show asymptotic bounds for $G_a$ and $G_b$ separately. 

Therefore, we only need to consider ``good" part where $\alpha_e = 1,2,3$. If $\alpha_e = 2$, the expression does not change if we substitute $\alpha_e$ with $\alpha_e - 2$. To see why, suppose $e = (a,b)$. Let us first discuss the effect of deleting such an edge. Firstly, the $D$'s terms in \eqref{e:tGcopy} will be multiplied by $D$ because one repeated ``good" edge is gone. Secondly, $w_{a}$ and $w_b$ remain the same. Thirdly, the length of the walk $(k)$ decreases by $2$, the term $D^{\frac{k}{2}}/n^{\frac{k}{2}}$ will also decrease by a factor of $n/D$. Thus, the overall bound will increase by a factor of at least $n$. 

In this case, one may worry if the graph will be disconnected into two connected components since we delete an edge. Denote the graph with $e$ deleted as $G'$. Note that $G'$ may be a connected graph or a disconnected graph, and either case works. The reason is best summarized with the following general decomposition lemma.
\begin{lemma}[Decomposition of graphs]
    \label{decomposition}
    Suppose $G = G_a \cup G_b$ is a $k$-step graph with vertex set $V$ and edge set $E$. $V_a$ (resp. $V_b$) is the vertex set of $G_a$ (resp. $G_b$). $E_a$ (resp. $E_b$) is the edge set of $G_a$ (resp. $G_b$). Suppose $E = E_a \sqcup E_b$ where $\sqcup$ denotes disjoint union. Suppose $V = V_a \cup V_b$. Then, the following holds
    $$t_G \leq t_{G_a} t_{G_b}.$$
\end{lemma}
\begin{proof}
    Observe the fact that $\sum_{l=1}^n d_l^{w_i} \leq \sum_{l=1}^n d_l^{w_i'}\sum_{l=1}^n d_l^{w_i''}$ where $w_i'+w_i''=w_i$ for all $i$. Then, we may write
    \begin{align*}
    t_{G}
    &= \frac{D^{\frac{k}{2}}}{n^{\frac{k}{2}}}
    \prod_{i \in V} \sum_{l = 1}^n d_l^{w_i}
    \frac{1}{D^rD^\frac{3p}{2}D^L}  \\
    &\leq \frac{D^{\frac{k_a}{2}}}{n^{\frac{k_a}{2}}}
    \prod_{i \in V_a} \sum_{l = 1}^n d_l^{w_i'}
    \frac{1}{D^{r_a}D^\frac{3p_a}{2}D^{L_a}}
    \times
    \frac{D^{\frac{k_b}{2}}}{n^{\frac{k_b}{2}}}
    \prod_{i \in V_b} \sum_{l = 1}^n d_l^{w_i''}
    \frac{1}{D^{r_b}D^\frac{3p_b}{2}D^{L_b}} = t_{G_a} t_{G_b}
\end{align*}
where $k_a$ (resp. $k_b$) is the number of steps in $G_a$ (resp. $G_b$) and $w_i'$ (resp. $w_i''$) is the corresponding term for $w$ in $G_a$ (resp. $G_b$). 
\end{proof}

\begin{remark}
    \label{remark-decomposition}
    We will mainly use the lemma in the following two scenarios: (1) when $G$ is a graph with one connected component but can be decomposed into two induced subgraphs $G_a$ and $G_b$; (2) when $G$ is a graph with two separate connected components $G_a$ and $G_b$. In both scenarios we have that $t_G \leq t_{G_a} t_{G_b}$ and it suffices to achieve bounds separately on $t_{G_a}$ and $t_{G_b}$. Note that the idea of decomposing a graph into several components or subgraphs will be frequently used in later proofs.
\end{remark}

% In the former case, we simply have $t_G \leq t_{G'}$. In the latter case, we denote $G' = G_a \cup G_b$. By the fact that $\sum_{l=1}^n d_l^{w_i} \leq \sum_{l=1}^n d_l^{w_i'}\sum_{l=1}^n d_l^{w_i''}$ where $w_i'+w_i''=w_i$ for all $i$, we can decompose $t_{G'}$ into $t_{G_a}t_{G_b}$. The reasoning is more explicit as below
% \begin{align*}
%     &\frac{D^{\frac{k}{2}}}{n^{\frac{k}{2}}}
%     \prod_{i \in V_k} \sum_{l = 1}^n d_l^{w_i}
%     \frac{1}{D^rD^\frac{3p}{2}D^L} = t_{G'} \\
%     &\leq \frac{D^{\frac{k_a}{2}}}{n^{\frac{k_a}{2}}}
%     \prod_{i \in V_k} \sum_{l = 1}^n d_l^{w_i'}
%     \frac{1}{D^{r_a}D^\frac{3p_a}{2}D^{L_a}}
%     \times
%     \frac{D^{\frac{k_b}{2}}}{n^{\frac{k_b}{2}}}
%     \prod_{i \in V_k} \sum_{l = 1}^n d_l^{w_i''}
%     \frac{1}{D^{r_b}D^\frac{3p_b}{2}D^{L_b}} = t_{G_a} t_{G_b}
% \end{align*}
% where we divide $k$ steps into $G_a$ and $G_b$. 
Now we go back to the claim about $G'$. If $G'$ is still a connected graph, we simply have $t_G \leq t_{G'}$. If $G'$ has two connected components $G_a$ and $G_b$, by Lemma \ref{decomposition} we know $t_{G} \leq t_{G'} \leq t_{G_a} t_{G_b}$. We can show asymptotic bounds for $G_a$ and $G_b$ separately. In the most special case where $(a,b)$ is the only edge that involves $b$, we can still delete it. The only difference is that now the term in \eqref{e:tGcopy} with $w_b = 1$, which has a contribution of $n$, will disappear, leading to a decrease by the factor of $1/n$. The net effect of reduction in this special case is thus $1$. Overall, we can freely drop the edges with multiplicities $2$. Note that we will not drop the entire graph because we are discussing graphs which have at least one odd multiplicity.

\subsubsection{Transformation of the ``bad" part}
% There may be several groups of ``bad" edges in the ``bad" part, but those can be dealt with separately. For simplicity we will address the case with one group of ``bad" edges here. 
The main idea here is ``disentangling" the ``bad" edges. For example, $(a,b)$ and $(b,c)$ are originally a pair of ``bad" edges, where both edges pick the same half-edge from $b$. To ``disentangle" $(a,b)$ and $(b,c)$ on $b$ means transforming to the same pair of edges $(a,b)$ and $(b,c)$ but with the difference that they are ``good" edges and that they pick two different half-edges from $b$. Note that we can also ``disentangle" edges with multiplicities larger than one. For example, $(a,b)$ with multiplicity $3$ and $(b,c)$ with multiplicity $2$ are ``bad" edges, sharing the same half-edge on $b$. In this case, only one half-edge from $b$ is used. To ``disentangle" them means that we consider the same pair of edges $(a,b)$ and $(b,c)$ with their corresponding multiplicities. However, two half-edges from $b$ are used, one with multiplicity $3$ for $(a,b)$ and the other with multiplicity $2$ for $(b,c)$.

Now we try to ``disentangle" the ``bad" edges as far as we can. At the same time we need to ensure that the bounds will not decrease during each step. Here we present two important observations on whether we can ``disentangle" the edges. 
\\\\
\textbf{Observation 1.} If the disentanglement does not produce any new ``good" edges, then we can conduct the disentanglement. \\
\textbf{Observation 2.} If all the ``good" edges produced by the disentanglement have multiplicities larger than 1, then we can conduct the disentanglement. 

\begin{proof}
    Observation 1 is not hard to verify. Suppose $(a,b)$ and $(b,c)$ are disentangled on $b$. After the disentanglement, $(a,b)$ and $(b,c)$ are still ``bad" edges sharing half-edges with some other edges. Since the disentanglement does not produce any new ``good" edges, then everything in the terms of $D$ in \eqref{e:tGcopy} remains the same. The only thing that changes is the power $w_b = z_b - \sum_{b \in e, e\in E_k} \alpha_e/2 + z_b' - \sum_{b \in e', e'\in E_k'}  \beta_{e'}/2$. $z_b'$ becomes larger because we are now picking one more distinct half-edge from $b$ in the ``bad" part. This change will only increase the bound thanks to Lemma \ref{lemmawi}.
    
    For Observation 2, suppose $(a,b)$ and $(b,c)$ with multiplicities $m_1$ and $m_2$ are disentangled on $b$. After the disentanglement, $(a,b)$ and $(b,c)$ both become ``good" edges with their multiplicities unchanged. Since the disentanglement changes two ``bad" edges into two ``good" repeated edges, the value of $L$ decreases by $2$ and the value of $r$ increases by $2$, which doesn't change the value of $D$'s terms in \eqref{e:tGcopy}. For $w$'s terms, again, the powers will only become larger because we are picking more distinct half-edges. Thus, this change will only increase the bound. 
    
    Note that Observation 2 also includes the case where only one ``good" edge is produced, which is easy to verify using the same logic. 
\end{proof}
The essence of the two observations is that any disentanglement will be valid as long as it does not produce any ``good" simple edges.

Now, we reduce the graphs by disentangling each ``bad" edge as far as we can, until any further disentanglement will violate the conditions in the two observations. The final transformed graph will result in only $3$ basic types of ``bad" sub-graphs. We call them $G1(m)$, $G2(m)$, and $G3(m)$ respectively. Note that all these 3 basic types are minimal, which is easier to verify. Also, any ``bad" sub-graph other than the 3 basics types can be further disentangled at some place.  
\\\\
\textbf{Type $\textbf{G1(m)}$.} ``Bad" graph of $(a,b)$ with multiplicity $1$ and $(b,c)$ with multiplicity $m$. They share a half-edge on $b$.\\
\textbf{Type $\textbf{G2(m)}$.} ``Bad" graph of $(a,b)$ with multiplicity $1$, $(b,c)$ with multiplicity $m$, and $(c,d) $ with multiplicity $1$. They share half-edges on $b$ and $c$. \\
\textbf{Type $\textbf{G3(m)}$.} ``Bad" graph of $(a,b_1)$ with multiplicity $1$, $(a,b_2)$ with multiplicity $1$, $\cdots$, $(a,b_m)$ with multiplicity $1$. They all share one half-edge on $a$. 
\\

Here, we prove that any ``bad" graph can be reduced into the union of $3$ basic types above. 
\begin{proposition}
    The ``bad" part of any minimal graph will only contain copies of the $3$ basic types above. A graph is said to be minimal if no further disentanglement can be performed on the graph.
\end{proposition}
\begin{proof}
    First, we show that the diameter of the minimal graph cannot exceed $3$. 
    
    For the purpose of contradiction, we assume that the diameter of the final graph is equal to or larger than $4$. By the definition of diameter, we know that there exist a chain of $(a,b), (b,c), (c,d),(d,e)$, which share half-edges on $b,c,d$. Thus, we can transform this final graph by disentangling $(b,c)$ and $(c,d)$ on $c$. This transformation will not produce any new ``good" edges. Thus, the original final graph is not minimal. Therefore, the diameter of the minimal graph cannot exceed $3$.

    Next, we consider all cases of graphs with diameter less than or equal to $3$ and check that they either fall into the $3$ basic types or are not minimal. 

    Consider a graph with diameter $2$. The most general form of such a graph will be a graph of $(a,b_1), (a,b_2), \cdots , (a,b_m)$ with multiplicities $x_1, x_2, \cdots , x_m$. They all share one half-edge on $a$. Firstly, assume that $m \geq 3$. If any $x_i$ is larger than $1$, we can disentangle it from the graph, because the disentanglement will only produce a repeated ``good" edge. By disentangling all repeated ``bad" edges, we will get exactly $G3(m')$. Next, assume that $m = 2$. If both are repeated ``bad" edges, then we can disentangle them because the disentanglement will only produce two repeated ``good" edges. If one is repeated and the other is simple, then if falls into the type $G1(m')$. If both are simple, then it falls into the type $G1(1)$. 

    Consider a graph with diameter $3$. If the graph contains a circle, we can of course disentangle on one vertex of the circle without producing new ``good" edges. Thus, we only need to consider acyclic graphs. The most general form will be a graph of $(a,b_1), (a,b_2), \cdots , (a,b_m)$ with multiplicities $x_1, x_2, \cdots , x_m$, $(a,c)$ with multiplicity $y$, and $(c,d_1), (c,d_2), \cdots , (c, d_n)$ with multiplicities $z_1, z_2, \cdots , z_n$. Note that here we abuse the notation a bit, but such notations will only be used in this proof. Now, for any $x_i$ and $z_i$ larger than $1$, we can directly disentangle those edges by the same logic. Without loss of generality, we consider $3$ overall cases: (1) $m\geq 2$ and $n\geq 2$ with all $x,z = 1$; (2) $m\geq 2$ and $n = 1$ with all $x = 1$ but $z_1$ being any integer; (3) $m = 1$ and $n = 1$ with $x_1,z_1$ being any integer. 

    In case (1), we first disentangle the graph on $a$ so that the left part becomes $G3$ and the right part becomes $(a,c)$ plus the cluster surrounding $c$. This disentanglement will not produce any ``good" edges. If $y = 1$, the remaining part is also $G3$ so we are done. If $y \geq 2$, we further disentangle $(a,c)$ from the cluster surrounding $c$ and the cluster becomes $G3$. This disentanglement will only produce a repeated ``good" edge $(a,c)$ so it's valid. 

    In case (2), we again disentangle the graph on $a$ so that the left part becomes $G3$ and the right part becomes $(a,c)$ plus $(c,d_1)$, both with general multiplicities. This disentanglement will not produce any ``good" edges. If $y\geq 2$ and $z_1\geq 2$, we can of course disentangle them on $c$. If at least one of $y$ and $z_1$ is $1$, it falls in type $G1$.

    In case (3), we simply have a chain of $(b_1,a)$, $(a,c)$, and $(c,d_1)$, each with arbitrary multiplicity. If all of them are repeated, we can of course disentangle all of them. If two of them are repeated, we further discuss their locations. If the two are adjacent, we disentangle on the vertex they share, producing one repeated ``good" edge and one $G1$. If the two are separate, we disentangle on $a$, again producing one repeated ``good" edge and one $G1$. If only one of them are repeated, we further discuss their locations. If the repeated edge is on the side ($(b_1,a)$ or $(c,d_1)$), we respectively disentangle on $a$ or $c$, producing one repeated ``good" edge and one $G1$. If the repeated edge is in the middle, it falls in type $G2$. If all of them are simple, it falls in type $G2$. 

    Therefore, we have proven that any minimal graph will fall in one of the $3$ basic types.
\end{proof}
Then we investigate the $3$ basic types and reduce the graph further. 

For Type $G1(m)$ where $m \geq 3$, we can always reduce $m$ by $2$. First of all, since $m \geq 3$, this reduction will not create any ``good" edges. Thus, the $D$'s terms in \eqref{e:tGcopy} remain the same. Also, $w_b$ and $w_c$ will both increase by $1$ because $\beta_{(b,c)}$ decreases by $2$. The middle term in \eqref{e:tGcopy} will thus be multiplied by at least $\sqrt{D/n} \times \sqrt{D/n}  = D/n$. Since the length of the walk $(k)$ decreases by $2$, the term $D^{\frac{k}{2}}/n^{\frac{k}{2}}$ will also decrease by a factor of $n/D$. Thus, the overall bound for $t_G$ in \eqref{e:tG} will either remain the same or increase. This reduction is valid.

For Type $G1(m)$ where $m = 2$, we can also reduce $m$ by $2$. An important difference is that the result will not be Type $G1$ any more. Instead, the reduced graph will be one ``good" edge. Firstly, the $D$'s terms will be multiplied by $D^2 \times D^{-\frac{3}{2}} = \sqrt{D}$ because two ``bad" edges are gone and one ``good" simple edge is present. Secondly, $w_b$ will increase by $1$ and $w_c$ will remain the same, which multiplies the product by at least $\sqrt{D/n}$. Thirdly, the length of the walk $(k)$ decreases by $2$, the term $D^{\frac{k}{2}}/n^{\frac{k}{2}}$ will also decrease by a factor of $n/D$. Thus, the overall bound  for $t_G$ in \eqref{e:tG}  will increase by a factor of at least $\sqrt{n}$. This reduction is valid.

For Type $G2(m)$ where $m \geq 3$, we can again reduce $m$ by $2$ freely. The reason is the same as Type $G1(m)$. 

For Type $G2(m)$ where $m = 2$, we can also reduce $m$ by $2$. An important difference is that the result will be two ``good" edges. Firstly, the $D$'s terms will be multiplied by $D^3 \times D^{-3} = 1$ because three ``bad" edges are gone and two ``good" simple edges are present. Secondly, $w_b$ and $w_c$ will both increase by $1$, which multiplies the product by at least $\sqrt{D/n} \times \sqrt{D/n} = D/n$. Thirdly, the length of the walk $(k)$ decreases by $2$, the term $D^{\frac{k}{2}}/n^{\frac{k}{2}}$ will also decrease by a factor of $n/D$. Thus, the overall bound  for $t_G$ in \eqref{e:tG}  will either remain the same or increase. This reduction is valid.

For Type $G3(m)$ where $m \geq 4$, we can reduce $m$ by $2$. We do this by removing $(a,b_{m})$ and $(a,b_{m-1})$ and by adding $(b_{m}, b_{m-1})$. Note that the reduced graph will be $G3(m-2)$ and a ``good" edge $(b_{m}, b_{m-1})$. Firstly, the $D$'s terms will be multiplied by $D^2 \times D^{-\frac{3}{2}} = \sqrt{D}$ because two ``bad" edges are gone and one ``good" simple edge is present. Secondly, $w_{a}$ will increase by $1$, which multiplies the product by at least $\sqrt{D/n}$. Thirdly, the length of the walk $(k)$ decreases by $1$, the term $D^{\frac{k}{2}}/n^{\frac{k}{2}}$ will also decrease by a factor of $\sqrt{n/D}$. Thus, the overall bound for $t_G$ in \eqref{e:tG} will increase by a factor of at least $\sqrt{D}$. This reduction is valid.

For Type $G3(3)$, we can also reduce $m$ by $2$. An important difference is that the result will be two ``good" edges $(a, b_1)$ and $(b_2, b_3)$. Firstly, the $D$'s terms will be multiplied by $D^3 \times D^{-3} = 1$ because three ``bad" edges are gone and two ``good" simple edges are present. Secondly, $w_{a}$ will increase by $1$, which multiplies the product by at least $\sqrt{D/n}$. Thirdly, the length of the walk $(k)$ decreases by $1$, the term $D^{\frac{k}{2}}/n^{\frac{k}{2}}$ will also decrease by a factor of $\sqrt{n/D}$. Thus, the overall bound for $t_G$ in \eqref{e:tG} will either remain the same or increase. This reduction is valid.

After the reductions above, the transformed graph's ``bad" part will only contain $G1(1)$ and $G2(1)$. 
\subsubsection{Contribution of the transformed graph}
Putting together everything, the current graph have $\alpha_e = 1, 3$ in the ``good" part and $G1(1), G2(1)$ in the ``bad" part. Since all degrees are even, and we have only been removing even numbers of edges, the graph can be decomposed into cycles and, since we allow parallel edges, walks of the form $a-b-a$. Again note that $\sum_{l=1}^n d_l^{w_i} \leq \sum_{l=1}^n d_l^{w_i'}\sum_{l=1}^n d_l^{w_i''}$ where $w_i'+w_i''=w_i$ for all $i$. Then, if we decompose $G$ into cycles $C_1,C_2,\cdots ,C_y$ and walks of the form $a-b-a$ denoted by $W_1, W_2,\cdots ,W_z$. We then have $t_G \leq t_{C_1}t_{C_2}\cdots t_{C_y}t_{W_1}t_{W_2}\cdots t_{W_z}$. From what we know before, any $t_{W_i}$ has an asymptotic contribution of $0$. For a cycle $C_i$, we have that $\alpha_e=1,3$ and $\beta_e = 1$. In any cycle, we have $w_i=1,-1,0$ since the two incoming edges in a cycle can have $\alpha_e=1,3$, so $\sum_{l=1}^n d_l^{w_i}\leq D$. 

Here is an important clarification. The $k$ we are considering does not correspond to the original $k$, but refers to the number of steps in the part we are investigating, such as the cycle $C_i$. As long as we show that each part (with their distinct $k$'s) goes to $0$ asymptotically, we can conclude that $t_G$ goes to $0$ asymptotically. 

Consider any circle $C$ that consists of ``good" simple edges, ``good" repeated edges with multiplicity $3$, $G1(1)$, and $G2(1)$. We need to show that such a circle will make an asymptotically small contribution to the expression. 

It is noteworthy that the four components of the circle $C$ are not equally costly.

If we change a $G2(1)$ into $G1(1)$, the bound for \eqref{e:tG} will only increase. Firstly, the $D$'s terms will be multiplied by $D$ because one ``bad" edge is gone. Secondly, one term with $w = 0$ will be gone because we deleted one vertex, which multiplies the product by $1/n$. Thirdly, the length of the walk $(k)$ decreases by $1$, the term $D^{\frac{k}{2}}/n^{\frac{k}{2}}$ will also decrease by a factor of $\sqrt{n/D}$. Thus, the overall bound will increase by a factor of at least $\sqrt{D/n}$. Therefore, we change all $G2(1)$ into $G1(1)$. 

If we change a ``good" simple edge into a ``good" repeated edge with multiplicity $3$, the bound will also only increase. Suppose the edge is $(a,b)$. Firstly, the $D$'s terms will be multiplied by $\sqrt{D}$ because one ``good" simple edge is gone and one ``good" repeated edge is present. Secondly, $w_a$ will decrease by $1$ and $w_b$ will decrease by $1$, which multiplies the product by $n/D$. Thirdly, the length of the walk $(k)$ increases by $2$, the term $D^{\frac{k}{2}}/n^{\frac{k}{2}}$ will also increase by a factor of $D/n$. Thus, the overall bound will increase by a factor of at least $\sqrt{D}$. Therefore, we change all ``good" simple edges into ``good" repeated edges with multiplicity $3$.

Now, we only have two components left. We will show that if the two consecutive units are of the same kind, we can combine them into one. 

Suppose that we have two consecutive ``good" repeated edges with multiplicities $3$, called $(a,b)$ and $(b,c)$ respectively. If we change them into $(a,c)$ with multiplicity $3$, the bound  \eqref{e:tG} will only increase. Firstly, the $D$'s terms will be multiplied by $D$ because one ``good" repeated edge is gone. Secondly, the term of $w_b = 2 - 6/2 = -1$ will disappear, which multiplies the product by $D^{\frac{1}{2}}n^{-\frac{3}{2}}$. Thirdly, the length of the walk $(k)$ decreases by $3$, the term $D^{\frac{k}{2}}/n^{\frac{k}{2}}$ will also decrease by a factor of $n^{\frac{3}{2}}D^{-\frac{3}{2}}$. Thus, the overall bound will remain the same. Therefore, we combine all consecutive ``good" repeated edges. 

Suppose that we have two consecutive $G1(1)$'s, called $(a, b)$, $(b,c)$, $(c,d)$, $(d,e)$ where $(a, b)$ and $(b,c)$ share a half-edge on $b$, $(c,d)$ and $(d,e)$ share a half-edge on $d$. If we change them into $(a,c)$ and $(c, e)$ which share a half-edge on $c$, the bound will only increase. Firstly, the $D$'s terms will be multiplied by $D^2$ because two ``bad" edges are gone. Secondly, the term of $w_b = 0$ and $w_d = 0$ will disappear, which multiplies the product by $1/n^2$. The term of $w_c$ will change from $1$ to $0$, which multiplies the product by $n/D$. Thirdly, the length of the walk $(k)$ decreases by $2$, the term $D^{\frac{k}{2}}/n^{\frac{k}{2}}$ will also decrease by a factor of $n/D$. Thus, the overall bound will remain the same. Therefore, we combine all consecutive $G1(1)$'s. 

Finally, we consider a cycle alternatively consisting of ``good" repeated edges with multiplicity $3$ and $G1(1)$. Suppose there are $a$ copies of ``good" repeated edges and $a$ copies of $G1(1)$. (Recall that the two components appear alternatively) Since what we consider is a $k$-walk, we have $3a + 2a = k$. Note that $k$ here may not correspond to the original $k$ since we have transformed the graph. However, we will show that any $k$ here will yield an asymptotic bound. Now consider the vertices within each $G1(1)$ and the vertices on the interface between ``good" repeated edge and $G1(1)$, it is not hard to see that all such vertices have $w = 0$. Finally, we can compute the bound for such a circle
\begin{align*}
    t_C = 
    \frac{D^{\frac{k}{2}}}{n^{\frac{k}{2}}}
    \prod_{i \in V_k} \sum_{l = 1}^n d_l^{w_i}
    \frac{1}{D^rD^\frac{3p}{2}D^L} 
    = \frac{D^{\frac{k}{2}}}{n^{\frac{k}{2}}} \frac{n^{2a + a}}{D^{2a + a}}
    = \left(\frac{n}{D}\right)^{3a - \frac{k}{2}}
    = \left(\frac{n}{D}\right)^{\frac{3k}{5} - \frac{k}{2}}
    = \left(\frac{n}{D}\right)^{\frac{k}{10} }
    \rightarrow 0.
\end{align*}
Thus, we conclude that all graphs $G$ which have at least an edge with odd multiplicity make contributions of $0$ asymptotically.  
% \\\\
% Abusing notation, if there are $r$ unique repeated edges (we call these $r$-edges) and $p$ unrepeated edges (we call these $p$-edges) in cycle $C_i$, we can show any cycle has a small contribution. To analyze cycles, first note that any edge must be an $r$-edge or a $p$-edge. Let the $r$-edges be distributed in the cycle in ``runs" of lengths $r_1,r_2,\cdots ,r_h$ such that $r_1+r_2+\cdots +r_h$. Between any two consecutive runs, there must be at least one $p$-edge. Since the vertices that contribute most are the vertices that have 2 $p$-edges, the largest contribution will occur when there is only one $p$-edge between runs, and all the other $p$-edges are distributed in a long run. Then, we have that
% \begin{align*}
%     \prod_{i \in V_k}\sum_{l=1}^nd_l^{w_i}&\leq\left(\sum_{l=1}^nd_l^{-1}\right)^{r-h}\left(\sum_{l=1}^nd_l^{0}\right)^{2h}\left(\sum_{l=1}^nd_l^{1}\right)^{p-h}\\
%     &\leq\left(\frac{n^\frac{3}{2}}{D^\frac{1}{2}}\right)^{r-h}n^{2h}D^{p-h}\\
%     &\leq\frac{n^\frac{3r}{2}}{D^\frac{r}{2}}D^p\frac{D^\frac{h}{2}}{n^\frac{3h}{2}}\frac{1}{D^h}\\
%     &= \frac{n^\frac{3r}{2}}{D^\frac{r}{2}}D^p\frac{1}{n^\frac{3h}{2}D^\frac{h}{2}}
% \end{align*}
% therefore cycles have the biggest contribution when $h$ is minimized, so if there are any $r$-edges, they should form one long run to maximize the contribution so cycles are maximized when $h=1$ (if $h=0$ there are only $p$-edges and it is simple to show that cycles are small in this case). Then, directly computing, the cycle's contribution is $\frac{1}{n^\frac{p+5}{2}D^\frac{1}{2}}$ which is small. 

\subsection{Conclusion of \eqref{e:moment_converge}}\label{subsec2.4}

% Note that any odd $k$ must fall into the above case, so that $\lim_{n\rightarrow \infty} \frac{1}{n}\E[\Tr M^k] = 0$. This proves one case of the result we are trying to show.

Suppose $k$ is odd. Then there must be an edge which has an odd multiplicity. As shown above, all such graphs have contribution of $0$ asymptotically. Thus, $\lim_{n\rightarrow \infty} \frac{1}{n}\E[\Tr M^k] = 0$ for all odd $k$. 

Suppose $k$ is even. As shown above, a graph with any edge which has an odd multiplicity has a contribution of $0$ asymptotically. Therefore, the only graphs that do not make small contributions must have even multiplicities for all edges. Among these graphs, we first investigate the ``good" graphs whose edges all have multiplicities $2$. Thus, $p=0$, $r=k/2$, and $L = 0$. Note that in this case, by Remark \ref{eqwhen2}, the equality is achieved for the bounds in Proposition \ref{goodrepeated}. Therefore, for any $G$ of this kind,
\begin{align*}
    &\frac{(1+o(1)) C_p C_L}{n} t_G  
    =  \frac{(1+o(1))}{n}\frac{D^{\frac{k}{2}}}{n^{\frac{k}{2}}} \prod_{i \in V_k}n\frac{1}{D^\frac{k}{2}}
    = (1+o(1))\frac{n^{|V_k|}}{n^{\frac{k}{2}+1}}.
\end{align*}
% \frac{(1+o(1)) C_p C_L}{n}\frac{D^{\frac{k}{2}}}{n^{\frac{k}{2}}}
%     \prod_{i \in V} \sum_{l = 1}^n d_l^{w_i}
%     \frac{1}{D^{r+\frac{3p}{2}+L}}
We now note that $|V_k| \leq |E_k|+1 = k/2+1$ with equality holding if and only if $G$ is a tree. If the equality does not hold, this term becomes asymptotically small. Therefore, the terms that make contributions to the sum are those that have $\alpha_e=2$ and are trees. Each of these makes a contribution of $1$. 

Then, let us consider the ``good" graphs whose edges not only have multiplicities $2$ but also larger even numbers. For those graphs, we may reduce the edges with larger even numbers by $2$ as we did before and get a reduced graph that have $\alpha_e = 2$ for all edges. The effect of reducing multiplicities is as follows. Suppose the reduced edge is $(a,b)$. Firstly, the length of the walk $(k)$ decreases by $2$ so the expression decreases by a factor of $n/D$. Secondly, both $w_a$ and $w_b$ increases by $1$, thus by Lemma \ref{lemmawi} the sum increases by a factor of $\omega (D/n)$. Thus, the overall term increases by a factor of $\omega (1)$. Since the reduced graphs have contributions of at most $1$, we know that the unreduced graphs which have multiplicities larger than $2$ must have contributions of $o(1)$. Thus, such graphs make a contribution of $0$ asymptotically. 

Finally, let us consider the ``bad" graphs with even multiplicities. Note that for all such graphs, we can always disentangle them and get a corresponding ``good" graph of the same shape. However, during the process of disentangling, the $w$'s terms will increase by at least $1$ because we are using more half-edges from the vertices. Thus, the bound will increase by a factor of at least $\sqrt{D/n}$. Since the disentangled graphs have contributions of at most $1$, we know that the ``bad" graphs must have contributions of $o(n/D)$. Thus, such graphs make a contribution of $0$ asymptotically.

To summarize, for even $k$, the only graphs that make contributions of $1$ are those rooted planar trees with lengths $k/2$. We conclude the proof for semicircle law. \nobreak\hfill$\qed$

\begin{remark} [sufficient and necessary condition for the convergence of moments]
    \label{rmk}
    The condition of ``$d_i \gg \sqrt{D/n}$ for all $i$" in Theorem \ref{semicircle} can be slightly relaxed to the weaker condition below $$\sum_l d_l^{w_i}  = o\left(n\sqrt{\frac{n}{D}}^{|w_i|}\right) \text{ for all } w_i \leq -1.$$
    We have proved in Lemma \ref{lemmawi} that the former condition implies the latter. Note that the latter cannot imply the former, where a counter example can be easily constructed by only having one $d_i = O(\sqrt{D/n})$. The latter condition is optimal in the sense that it is the sufficient and necessary condition for the semicircle law to hold.

    For the sufficient part, it is not hard to see that during the entire proof, we only used the latter condition. The only part we used the former condition is in Lemma \ref{lemmawi} where we showed that it implies the latter condition. Therefore, the latter condition is sufficient for the semicircle law.

    For the necessary part, we prove it below. Assume that there is a $w_i \leq -1$ such that $\sum_l d_l^{w_i}  = \Omega(n({n/D})^{|w_i|/2})$. Denote $k = 2(1-w_i)$. Now consider a graph consisting of two vertices and one edge between them with multiplicity $k$. The contribution of such a graph is
    \begin{align*}
        &\frac{D^{\frac{k}{2}}}{n^{1+\frac{k}{2}}} (\sum_l d_l^{1-\frac{k}{2}})^2 \frac{1}{D}
        = \frac{D^{\frac{k}{2}}}{n^{1+\frac{k}{2}}} \Omega(n^2 \frac{n^{\frac{k}{2}-1}}{D^{\frac{k}{2}-1}}) \frac{1}{D}
        = \Omega(1).
    \end{align*}
    Thus, graphs of the above shape will also make contributions of $\Omega(1)$ even if they are not graphs with all multiplicities $2$. In other words, semicircle law does not hold without the optimal condition. Therefore, the latter condition is necessary for the semicircle law to hold.
\end{remark}

\subsection{Proof of Theorem \ref{concentration}} \label{subsec2.5}

In this subsection, we prove Theorem \ref{concentration}, which basically states that the moments of the empirical eigenvalue distribution concentrate around their expectations. To show it, we prove that 
\begin{align*}
    \E \biggl[ \Big|\frac{1}{n}\Tr[M^k] - \frac{1}{n}\E[\Tr[M^k]] \Big|^2  \biggr] \to 0.
\end{align*}
Moreover, the above quantity is the same as the variance, so it suffices to prove that
\begin{align*}
    \frac{1}{n^2}\left(\E\left[ \Tr[M^k]^2\right] - \E[\Tr[M^k]]^2\right) \to 0.
\end{align*}
Now we denote the term $(H_{i_ji_k}^{s_lt_m} - \frac{1}{D-1})$ as $h_{i_ji_k}^{s_lt_m}$. By the same computations, we have that
\begin{align*}
    & \frac{1}{n^2}\left(\E\left[ \Tr[M^k]^2\right] - \E[\Tr[M^k]]^2\right)\\
    &= \frac{D^k}{n^{k+2}} \biggl( 
    \sum_{\substack{i_1, \cdots , i_k, i_1', \cdots , i_k'\\s_1, t_2, \cdots , s_k, t_1, s_1', t_2',\cdots ,s_k', t_1'}} \frac{1}{d_{i_1}\cdots d_{i_k}d_{i_1'}\cdots d_{i_k'}} \E \biggl[ h_{i_1i_2}^{s_1t_2} \cdots  h_{i_ki_1}^{s_kt_1} h_{i_1'i_2'}^{s_1't_2'} \cdots  h_{i_k'i_1'}^{s_k't_1'} 
    \biggr]\\
    &- \sum_{\substack{i_1, \cdots , i_k, i_1', \cdots , i_k'\\s_1, t_2, \cdots , s_k, t_1, s_1', t_2',\cdots ,s_k', t_1'}} \frac{1}{d_{i_1}\cdots d_{i_k}d_{i_1'}\cdots d_{i_k'}} \E \biggl[ h_{i_1i_2}^{s_1t_2} \cdots  h_{i_ki_1}^{s_kt_1}\biggr] \E\biggl[ h_{i_1'i_2'}^{s_1't_2'} \cdots  h_{i_k'i_1'}^{s_k't_1'} \biggr]
    \biggr).
\end{align*}
Let $\mathcal{G}_{II'}$ denote the graph associated with the term $\E [ h_{i_1i_2}^{s_1t_2} \cdots  h_{i_ki_1}^{s_kt_1} h_{i_1'i_2'}^{s_1't_2'} \cdots  h_{i_k'i_1'}^{s_k't_1'} ]$, $\mathcal{G}_{I}$ associated with $\E [ h_{i_1i_2}^{s_1t_2} \cdots  h_{i_ki_1}^{s_kt_1}]$ and $\mathcal{G}_{I'}$ associated with $\E [ h_{i_1'i_2'}^{s_1't_2'} \cdots  h_{i_k'i_1'}^{s_k't_1'} ]$. 
% TODO: Add justification because this is two trees
From the analysis before, the maximum of the second term is achieved when $\mathcal{G}_{I}$ and $\mathcal{G}_{I'}$ both have $k/2 + 1$ vertices and $k/2$ ``good" edges of multiplicities $2$, without sharing any common vertex.

For the first term, we see that $\mathcal{G}_{II'}$ has at most $2$ connected components. However, by Lemma \ref{decomposition} and Remark \ref{remark-decomposition}, we can simply do the same analysis to the $2$ connected components. By repeating all the transformation and reduction as in the previous proof, we note that the only contributing graphs are when all edges have multiplicities of $2$. In this scenario, $\mathcal{G}_{II'}$ has at most $k$ distinct edges. Accompanied with the fact that $\mathcal{G}_{II'}$ has at most $2$ connected components, it has at most $k + 2$ vertices. Therefore the terms are not asymptotically vanishing if and only if $\mathcal{G}_{II'} = \mathcal{G}_{I} + \mathcal{G}_{I'}$ where $\mathcal{G}_{I} \cap \mathcal{G}_{I'} = \emptyset$ and they both have $k/2 + 1$ vertices and $k/2$ edges. This implies that the leading terms of the two terms in the expression above cancels, and the contribution from remaining terms is $o(1)$, which concludes the proof.

% we see that $\mathcal{G}_{II'}$ has at most $k$ edges and at most $2$ connected components. Therefore, $\mathcal{G}_{II'}$ has at most $k + 2$ vertices and the equality holds if and only if $\mathcal{G}_{II'} = \mathcal{G}_{I} + \mathcal{G}_{I'}$ and $\mathcal{G}_{I} \cap \mathcal{G}_{I'} = \emptyset$. This implies that the leading terms of the two terms in the expression above cancels, and the contribution from remaining terms is $o(1)$, which concludes the proof.

\section{Proof of Lemma \ref{equivlemma} and \Cref{equiv}} \label{sec:equivlemma}

\subsection{Proof of Lemma \ref{equivlemma}}

We generate a random graph $\cG$ with degree sequence ${ d_i }_{1 \leq i \leq n}$ using the configuration model as follows:
\begin{enumerate}
\item Begin with $n$ vertices, where each vertex $i$ has $d_i$ half-edges. Index these half-edges sequentially from $1$ to $D$.
\item While there exist vertices with fewer than $C\sqrt{D/n}$ unused half-edges, choose the unused half-edge with the smallest index from these vertices (if there exists any) and match it uniformly at random with another unused half-edge.
\item Repeat step 2 until every vertex has either $0$ or at least $C\sqrt{D/n}$ unused half-edges.
\item Let $\{ \tilde d_i \}_{1 \leq i \leq n}$ denote the number of unused half-edges at each vertex, so that either $\tilde d_i = 0$ or $\tilde d_i \geq C\sqrt{D/n}$ for each $i$. Matching these remaining half-edges yields a graph $\tilde \cG$, which is a configuration model graph with degree sequence $\{ \tilde d_i \}_{1 \leq i \leq n}$.
\end{enumerate}

In the above procedure, the graph $\tilde \cG$ can be obtained from $\cG$ by sequentially removing edges. In the following, we show the graph $\tilde \cG$ satisfies the properties in Lemma \ref{equivlemma}. We remark that 
condition on the degree sequence $\{ \tilde d_i \}_{1 \leq i \leq n}$, the law of 
$\tilde \cG$ is the same as a configuration model graph with that degree sequence. And the third statement of 
Lemma \ref{equivlemma} holds.

 To get the degree sequence $\{ \tilde{d}_i \}_{1\leq i\leq n}$ (viewed as $n$ vertices with each vertex $i$ has $\tilde d_i$ half-edges) from the degree sequence $\{ d_i \}_{1\leq i\leq n}$ in the procedure above, 
we have the stochastic process: 
choose the half-edge with the smallest index from vertices with fewer than $C\sqrt{D/n}$ half-edges, remove it, and another uniform randomly picked half-edge.  We continue this process of dropping half-edges in pairs of two until all vertices have degrees either $0$ or at least $C \sqrt{D/n}$. 
% we drop all vertices with degrees less than $C \sqrt{D/n}$ by eliminating all half-edges coming out of them. By doing this, there may be new vertices whose degrees decrease below $C \sqrt{D/n}$. We drop these vertices as well. We continue the process until all vertices have degrees larger than $C \sqrt{D/n}$. In the end, there will be a pruned graph $\tilde{\mathcal{G}}$ whose degrees are all larger than $C \sqrt{D/n}$. Suppose the pruned degree sequence is $\{ \tilde{d}_n \}$. An immediate fact is that, conditioned on $\{ \tilde{d}_n \}$, the law of the pruned graph $\tilde{\mathcal{G}}$ is still the configuration model. 

For the conciseness of notation, we call vertices with degrees $0$ or larger than $C\sqrt{D/n}$ \textbf{heavy vertices} and other vertices \textbf{light vertices}. We call half edges associated with heavy vertices \textbf{heavy half-edges}, and half edges associated with light vertices \textbf{light half-edges}.

Note that the difficulty of this task lies in the fact that the action of eliminating light vertices will also eliminate a fraction of half-edges connected to heavy vertices and transfer them into light vertices. We want to make sure that this pruning process will not lead to the situation that each vertex has degree $0$. Additionally, we specify the relative magnitudes of some asymptotic terms
\begin{equation}
    \label{rule}
    \lim_{n\to \infty} \varepsilon_n = 0,
\quad
    \lim_{n\to \infty} C_n = \infty,
\quad
    \lim_{n\to \infty} K_n = \infty,
\quad
    \lim_{n\to \infty} K_n\varepsilon_n C_n^2 = 0.
\end{equation}
Although the terms above have asymptotic behaviors, in this proof we will deal with a fixed large enough $n$ and thus fixed $\varepsilon_n, C_n, K_n$. Therefore, we simply write $\varepsilon_n, C_n, K_n$ as $\varepsilon, C, K$. 

For the rigorousness of proof, we formally model the pruning process. The pruning procedure is essentially removing half-edges in pairs of two, one from a light vertex and another one uniformly out of all half-edges. We repeat this step at each integer time stamp $t\in \mathbb{N}$ until there are no light vertices. The number of light half-edges at each time $t$: $\{ S_t\} _{t \in \mathbb{N}}$, can be viewed as a birth-death process. Particularly, $S_0 > 0$ is the number of light half-edges of the original graph. We write $S_t = S_0 + \sum_{i=1}^t X_i$ where $X_i$'s are the increments or decrements of light half-edges at each time $i$. Since the process stops when there are no light vertices, we use $\tau$ to denote the stopping time 
\begin{align}\label{e:deftau}
    \tau = \inf \{ t\in \mathbb{N}: S_t = 0\} \wedge 3\varepsilon C \sqrt{nD}.
\end{align}
The second part in \eqref{e:deftau} is to make sure $\tau\leq 3\varepsilon C \sqrt{nD}$. We will show with high probability $\tau\leq 3\varepsilon C \sqrt{nD}$.
Note that $\tau$ is half of the total number of deleted half-edges because at each $t$ we delete two half-edges. 

% å äžæ­¥\tau's stopping condition: \tau = 2\ep n
% Show that \tau != 2\ep n and thus S_\tau = 0

Now we consider the distribution of $X_t$, which is the net change of $S_t$ when we delete one half-edge from a light vertex and another one uniformly out of all half-edges. We denote $D_t$ the number of vertices which have degrees of $C\sqrt{D/n}$ at time $t$.
Firstly, we delete one half-edge from a light vertex. It implies that $S_t$ will decrease by $1$. Secondly, we randomly delete a half-edge, which comes from a vertex $v$. There are three cases:
\begin{enumerate}
  \item If $v$ happens to be a light vertex, then we know that $X_t = -1 - 1 = -2$. 
  \item If $v$ happens to be a vertex with degree $C\sqrt{D/n}$, then we will have to add $C\sqrt{D/n} - 1$ to $S_t$ because $v$ has transferred from a heavy vertex to a light vertex. In this case, $X_t = -1 + C\sqrt{D/n} - 1 = C\sqrt{D/n} - 2$. 
    \item If $v$ does not fit into the last two cases above, then $X_t = -1$. 
    \end{enumerate}
    In summary, 
$$X_t = 
\begin{cases} 
    -2 & \text{with probability } \frac{S_t}{\cor D-2t+1}, \\
    % & \\
    C\sqrt{\frac{D}{n}} -2 & \text{with probability } \frac{D_t \times C\sqrt{\frac{D}{n}}}{\cor D-2t+1}, \\
    -1 & \text{with probability } 1 - \frac{S_t}{\cor D-2t+1} - \frac{D_t \times C\sqrt{\frac{D}{n}}}{\cor D-2t+1}, \\
\end{cases}$$
where $D_t$ is the number of vertices which have degrees of $C\sqrt{D/n}$ at time $t$. 

Next we show that almost surely, for $t\leq \tau$, the random variable $D_t\leq 2 \varepsilon n$. The meaning of the upper bound is that at any time $t$, the number of vertices which have degrees of $C\sqrt{D/n}$ will be upper bounded by $2 \varepsilon n$. The reason is as follows. Recall our definition of $\tau$ from \eqref{e:deftau}, we have $\tau \leq 3\varepsilon C \sqrt{nD}$. Thus we are guaranteed that the total number of heavy half-edges deleted in this process will not exceed $3\varepsilon C \sqrt{nD}$. We divide all vertices into two groups, those with degrees smaller than $4C\sqrt{D/n}$ and those with degrees larger than $4C\sqrt{D/n}$. For the first group, we conservatively assume that all degrees are exactly $C\sqrt{D/n}$ at time $t$. The cardinality of the first group can be bounded by $\varepsilon n$ by simply letting $C = 4C$ in our assumption, provided $n$ is large enough. For the second group, if any vertex has its degree decreased to $C\sqrt{D/n}$, a number of at least $3C\sqrt{D/n}$ heavy half-edges must be deleted from it. Since the total number of heavy half-edges deleted will not exceed $3\varepsilon C \sqrt{nD}$, we know that the number of vertices which have degrees of $C\sqrt{D/n}$ that are reduced from the second group will not exceed $\frac{3\varepsilon C \sqrt{nD}}{3C\sqrt{D/n}} = \varepsilon n$. Thus, the total number of vertices which have degrees of $C\sqrt{D/n}$ will be upper bounded by $\varepsilon n + \varepsilon n = 2 \varepsilon n$ at any time $t$.

We can construct a new random variable $\tilde{X}_t$
$$\tilde{X}_t = 
\begin{cases} 
    C\sqrt{\frac{D}{n}} -1 & \text{with probability } \frac{4 \varepsilon n \times C\sqrt{\frac{D}{n}}}{D}, \\
    -1 & \text{with probability } 1 -\frac{4 \varepsilon n \times C\sqrt{\frac{D}{n}}}{D}, \\
\end{cases}$$
Then from the discussion above, almost surely, we have 
\begin{align*}
    \bP(\tilde X_t=C\sqrt{D/n}-1)=\frac{4 \varepsilon n \times C\sqrt{\frac{D}{n}}}{D}
    \geq \frac{2D_t \times C\sqrt{\frac{D}{n}}}{D}\geq \frac{D_t \times C\sqrt{\frac{D}{n}}}{\cor D-2t+1}= \bP( X_t=C\sqrt{D/n}-1).
\end{align*}
We can couple $\tilde X_t$ and $X_t$, such that $\tilde X_t\geq X_t$ almost surely. 
We define $\tilde{S}_t = \tilde{S}_0 + \sum_{i=1}^t \tilde{X}_i$. It immediately follows that $S_t \leq \tilde{S}_t$ almost surely for all $t$. We define the same stopping time 
\begin{align}\label{e:deftautilde}
    \tilde{\tau} = \inf \{ t\in \mathbb{N}: \tilde{S}_t = 0 \} \wedge 3\varepsilon C \sqrt{nD},
\end{align}
where 
\begin{align}\label{e:taucompare}
    \tilde{\tau} \geq \tau
\end{align} almost surely. Therefore, as long as we bound the stopping time $\tilde{\tau}$, we obtain an upper bound on the total number of deleted half-edges in the original process. 

% Note that we will finally let $\varepsilon$ go to $0$ and $C$ go to $\infty$. Since our choice of $\varepsilon$ and $C$ is arbitrary, we may choose their relative rates for the purpose of our proof. We assume that 
% $\varepsilon C^2 = \frac{1}{8}$ and thus that $\varepsilon^2 C^2$ will go to $0$.

Now, let us investigate the behavior of the process $\{ \tilde{S}_t\} _{t \in \mathbb{N}}$. Considering the construction above, we know that $\tilde{X}_t$ will be $-1$ for most of the time and an upsurge occasionally, i.e. $\tilde X_t=C\sqrt{D/n}-1$. We treat the time interval $t \in [1,  T_1]$ with $T_1=\varepsilon C \sqrt{nD}$ as the first round, and define the number of upsurges $\tilde{P}_1$ inside first round:
\begin{align}
    \tilde P_1:=\{1\leq t\leq T_1: \tilde X_t=C\sqrt{D/n}-1\}.
\end{align}
Then after the first round, we have
\begin{align}
    \tilde S_{T_1}=\tilde S_0+\tilde{P}_1 (C\sqrt{D/n}-1)-(T_1-\tilde P_1)=\tilde{P}_1 C\sqrt{D/n},
\end{align}
which is the initial value for the second round, $\tilde{P}_1 C\sqrt{D/n}$. It is noteworthy that the distribution of $\tilde{P}_1$ will be $Bin(\varepsilon C \sqrt{nD}, \frac{4 \varepsilon n  C\sqrt{D/n}}{D}) = Bin (\varepsilon C \sqrt{nD}, 4\varepsilon C \sqrt{n/D})$. 

Then we treat the time interval $[T_1+1, T_2]$ with $T_2-T_1=\tilde{P}_1 C\sqrt{D/n}$ as the second round.
Similarly, we define the number of upsurges during the second round as $\tilde{P}_2$ and know that 
$S_{T_2}=\tilde{P}_2 C\sqrt{D/n}$, and 
$\tilde{P}_2 \sim Bin( \tilde{P}_1 C\sqrt{D/n}, \frac{4 \varepsilon n  C\sqrt{D/n}}{D}) $. We repeat this process. We treat the time interval $[T_{k-1}+1, T_k]$ with $T_k-T_{k-1}=\tilde{P}_k C\sqrt{D/n}$ as the $k$-round, and  define the number of upsurges during the $k$-round as $\tilde{P}_k$
\begin{align}
    \tilde P_k:=\{T_{k-1}+1\leq t\leq T_k: \tilde X_t=C\sqrt{D/n}-1\}.
\end{align}
Then
$S_{T_k}=\tilde{P}_k C\sqrt{D/n}$, and 
$\tilde{P}_k \sim Bin( \tilde{P}_k C\sqrt{D/n}, \frac{4 \varepsilon n  C\sqrt{D/n}}{D}) $.

%{\color{red}until the $N$-th round, for $t\in [T_{N-1}+1, T_N]$, where there will be no upsurge with high probabilities. Notationally, we stop when $\p(\tilde{P}_N \geq 1 )$ is small for some $N$. That is also the stopping time for the process $\tilde{S}_t$.}

Now that preparations are set up, we start to discuss the probability statements regarding the process. In the first round, we sample $\tilde{P}_1 \sim Bin (\varepsilon C \sqrt{nD}, 4\varepsilon C \sqrt{n/D})$ and decide the starting point for the second round. By properties of the Binomial distribution and Chernoff bounds, we may obtain the expectation and the tail bound as follows,
$$\E \tilde{P}_1 = \varepsilon C \sqrt{nD} \times 4\varepsilon C \sqrt{\frac{n}{D}} = 4 \varepsilon^2 C^2 n =: \mu_1,$$
and it follows
$$\p (\tilde{P}_1 > 2\mu_1) \leq \exp (- \frac{1}{3} 4 \varepsilon^2 C^2 n).$$

For the second round, we condition on $\tilde{P}_1 \leq  2\mu_1 = 8 \varepsilon^2 C^2 n$. Then we sample the Binomial distribution $\tilde{P}_2 \sim Bin(\tilde{P}_1 C\sqrt{D/n},4\varepsilon C \sqrt{n/D}) $ where we know that $\tilde{P}_1 \leq 8 \varepsilon^2 C^2 n$. 
Let $\hat{P}_2 \sim Bin (8\varepsilon^2 C^3 \sqrt{nD}, 4\varepsilon C \sqrt{n/D})$, and we can couple it with $\tilde P_2$ such that $\tilde P_2\leq \hat P_2$ almost surely.
Now we may obtain the expectation the the tail bound as follows,
\begin{align}
    &\E [\tilde{P}_2|\tilde{P}_1 \leq 2\mu_1] \leq \E \hat{P}_2 = 8\varepsilon^2 C^3 \sqrt{nD} \times 4\varepsilon C \sqrt{\frac{n}{D}} = 32 \varepsilon^3 C^4 n =: \mu_2 \\
&\p(\tilde{P}_2 > 2\mu_2|\tilde{P}_1 \leq 2\mu_1) \leq \p(\hat{P_2} > 2\E \hat{P}_2) \leq \exp (-\frac{1}{3} 32 \varepsilon^3 C^4 n)
\end{align}

We repeat this process. In the $k$-th round, we condition on $\tilde{P}_{k-1} \leq 2\mu_{k-1} = (8 \varepsilon C^2)^{k-1} \varepsilon n$. Then we sample $\tilde{P}_k \sim Bin (\tilde{P}_{k-1} C\sqrt{D/n} , 4\varepsilon C \sqrt{n/D})$ where we know that $\tilde{P}_{k-1} \leq (8 \varepsilon C^2)^{k-1} \varepsilon n$. We couple $\tilde{P}_k$ with $\hat{P}_k \sim Bin ((8 \varepsilon C^2)^{k-1} \varepsilon C \sqrt{nD}, 4\varepsilon C \sqrt{n/D})$, such that $\hat{P}_k \geq \tilde{P}_k$ almost surely. This gives us the following  expectation, and the tail bound,
\begin{align}\begin{split}\label{e:defuk}
    &\E [\tilde{P}_k|\tilde{P}_{k-1} \leq 2\mu_{k-1}] \leq \E \hat{P}_k = (8 \varepsilon C^2)^{k-1} \varepsilon C \sqrt{nD} \times 4\varepsilon C \sqrt{\frac{n}{D}} = 2^{k-1} (4\varepsilon C^2)^k \varepsilon n =: \mu_k, \\
&\p(\tilde{P}_k > 2\mu_k |\tilde{P}_{k-1} \leq 2\mu_{k-1}) \leq \p(\hat{P}_k > 2\E \hat{P}_k) \leq \exp (-\frac{1}{3} 2^{k-1} (4\varepsilon C^2)^k \varepsilon n ).
\end{split}\end{align}

We consider the $N$-th round such that $2\mu_{N-1} = K$ where $K$ is a large number whose magnitude is specified in \eqref{rule}. Recall that $N$-th round is for time $t\in [T_{N-1}+1, T_N]$, where 
\begin{align}\label{e:recall_TN}
    T_N=(\varepsilon n+ \tilde P_1+\tilde P_2+\cdots+\tilde P_{N-1})\times C\sqrt{D/n}.
\end{align}
In the $N$-th round, we condition on $\tilde{P}_{N-1} \leq 2\mu_{N-1} = K$. Then we sample the binomial random variable $\tilde{P}_N \sim Bin (\tilde{P}_{N-1} C\sqrt{D/n} , 4\varepsilon C \sqrt{n/D})$ where we know that $\tilde{P}_{N-1} \leq K$. Again we can couple $\tilde{P}_N $ with $\hat{P}_N \sim Bin (K C\sqrt{D/n}, 4\varepsilon C \sqrt{n/D})$ such that $\hat{P}_N \geq \tilde{P}_N$ almost surely. This gives the following probability
\begin{align}\begin{split}\label{e:bound_PN}
    \p (\tilde{P}_N \geq 1 |\tilde{P}_{N-1} \leq 2\mu_{N-1} ) 
    \leq \p(\hat{P}_N \geq 1) 
    &= 1 - (1 - 4\varepsilon C \sqrt{\frac{n}{D}})^{KC\sqrt{\frac{D}{n}}} \\
    &= 1 - e^{- 4K\varepsilon C^2}
    \leq 4K\varepsilon C^2 .
\end{split}\end{align}
Here the second last step follows from the definition of $e$ and the fact that $\lim_{n\to \infty}\sqrt{D/n} = \infty $; and the last step follows from \eqref{rule} and a first-order Taylor expansion. 

Recall $\mu_k$ from \eqref{e:defuk}. We compute the following quantity, which will be useful later 
\begin{align}\begin{split}
    \varepsilon n + 2\mu_1 + 2\mu_2 + \cdots + 2\mu_{N-1} + 1 
    \leq &\varepsilon n + 2\mu_1 + 2\mu_2 + \cdots + 2\mu_{N-1} + 2\mu_{N} + \cdots\\
    =& \varepsilon n + 8\varepsilon^2 C^2 n + 64\varepsilon^3 C^4 n + \cdots \\
    \leq& \frac{1}{1 - 8\varepsilon C^2} \varepsilon n
    \leq 2\varepsilon n \label{e:sumMu}
\end{split}
\end{align}
where we know that $\varepsilon C^2$ will go to $0$ by using \eqref{rule}. Now we use $\tilde{P}=\tilde P_1+\tilde P_2+\cdots+\tilde P_N$ to denote the total number of upsurges. Then it satisfies
\begin{align}\begin{split}\label{e:tPbound}
    &\p (\tilde{P} > 2(\mu_1 + \mu_2 + \cdots + \mu_{N-1})+1 ) \\
    \leq &\p(\tilde{P}_1 > 2\mu_1) + \p(\tilde{P}_1 \leq 2\mu_1, \tilde{P}_2 > 2\mu_2) + \cdots + \p(\tilde{P}_1 \leq 2\mu_1, \tilde{P}_2 \leq 2\mu_2, \cdots, \tilde{P}_N \geq 1)\\
    \leq& \p(\tilde{P}_1 > 2\mu_1) + \p(\tilde{P}_2 > 2\mu_2|\tilde{P}_1 \leq 2\mu_1) + \cdots + \p (\tilde{P}_N \geq 1|\tilde{P}_{N-1} \leq 2\mu_{N-1}) \\
    \leq& \exp (- \frac{1}{3} 4 \varepsilon^2 C^2 n) + \exp (-\frac{1}{3} 32 \varepsilon^3 C^4 n) + \cdots + \exp (-\frac{1}{3} 2^{N-2} (4\varepsilon C^2)^{N-1} \varepsilon n ) +  4K\varepsilon C^2 \\
    =& \exp (- \frac{1}{3} 4 \varepsilon^2 C^2 n) + \exp (-\frac{1}{3} 32 \varepsilon^3 C^4 n) + \cdots + \exp (-\frac{1}{3} K ) +  4K\varepsilon C^2 \\
    \leq& \frac{1}{1 - 8\varepsilon C^2} \exp (-\frac{1}{3}K) +  4K\varepsilon C^2 \\
    \leq& 2\exp (-\frac{1}{3}K) +  4K\varepsilon C^2 ,
\end{split}\end{align}
where we used \eqref{e:defuk} and \eqref{e:bound_PN} for the fourth line.

In the following we denote the event 
\begin{align}
    \Omega:=\{\tilde{P} \leq  2(\mu_1 + \mu_2 + \cdots + \mu_{N-1})+1 \}.
\end{align}
Then \eqref{e:tPbound} implies
\begin{align}
   \p(\Omega)\geq 1-2\exp (-\frac{1}{3}K) -  4K\varepsilon C^2=1-o(1).
\end{align}
Moreover, on $\Omega$, we have
\begin{align}
    S_{T_N}=\tilde P_NC\sqrt{D/n}=0.
\end{align}
Recall that $\tilde{\tau}$ from \eqref{e:deftautilde} is the stopping time for the process, then $\tilde \tau\leq T_N$. Recall $T_N$ from \eqref{e:recall_TN}, it follows 
\begin{align}
\begin{split}
    \tilde \tau&\leq T_N\leq (\varepsilon n+ \tilde P)\times C\sqrt{D/n}
    \leq (\varepsilon n + 2(\mu_1 + \mu_2 + \cdots + \mu_{N-1})+1)\times C\sqrt{D/n}\\
    &\leq 2\varepsilon n \times C\sqrt{D/n}
\end{split}
\end{align}
where the last inequality follows from \eqref{e:sumMu}.
By using the fact that $\tilde{\tau} \geq \tau$ almost surely (recall from \eqref{e:taucompare}), we conclude that on $\Omega$, 
\begin{align*}
   \tau \leq \tilde \tau\leq 2\varepsilon  \times C\sqrt{nD}.
\end{align*}
Therefore, we conclude that the total number of half-edges deleted will be less than $4\varepsilon C \sqrt{nD}$ condition on $\Omega$, because $\tau$ is half of the total number of half-edges deleted. Meanwhile, the number of vertices deleted will be less than $ \varepsilon n+\tilde P\leq  2\varepsilon n$ condition on $\Omega$.  This finishes the proof of \Cref{equivlemma}.

\begin{comment}
Therefore, we have $\tilde{\tau} = (\varepsilon n + \tilde{P}) \times C\sqrt{D/n}$. Therefore, 
\begin{align*}
    &\p (\tilde{\tau} > (\varepsilon n + 2(\mu_1 + \mu_2 + \cdots + \mu_{N-1}))\times C\sqrt{\frac{D}{n}}) \\
    = &\p(\tilde{P} > 2(\mu_1 + \mu_2 + \cdots + \mu_{N-1})) 
    \leq 2\exp (-\frac{1}{3}K) +  4K\varepsilon C^2 .
\end{align*}
Using the inequality \eqref{e:sumMu}, we know that
\begin{align*}
    \p (\tilde{\tau} > 2\varepsilon n \times C\sqrt{\frac{D}{n}}) 
    \leq 2\exp (-\frac{1}{3}K) +  4K\varepsilon C^2 .
\end{align*}

% 
Note that the RHS is asymptotically vanishing because of \eqref{rule}. Therefore, we conclude that the total number of half-edges deleted will be less than $4\varepsilon C \sqrt{nD}$ with probability $1 - o(1)$ because $\tau$ is half of the total number of half-edges deleted. Meanwhile, although we have two stopping conditions, it indicates that our process stops because $S_\tau = 0$ instead of $\tau = 3\varepsilon C \sqrt{nD}$. Similarly, the number of vertices deleted will be less than $\varepsilon n + 2\varepsilon n = 3\varepsilon n$ with probability $1-o(1)$. 
\end{comment}

\subsection{Proof of \Cref{equiv}}
% Since the proportion of the degrees less than $C\sqrt{\frac{D^{(n)}}{n}}$ is only $\varepsilon$, we drop this proportion of the degrees. 
In this section, we define the following notation for matrices: given any matrix $ A $ and any index set $ \mathcal{I} $, let $ A^{\mathcal{I}} $ denote the matrix obtained by removing the columns and rows of $ A $ indexed by $ \mathcal{I} $.

By Lemma \ref{equivlemma}, we know that the original graph $\mathcal{G}$ may be pruned using the graph pruning process and transferred into a graph $\tilde{\mathcal{G}}$ where all degrees are now larger than $C\sqrt{D/n}$. It also follows that, with probability $1 - o(1)$, the number of removed vertices is less than $2\varepsilon n$ and the number of removed half-edges is less than $4\varepsilon C \sqrt{nD}$. The law for the pruned graph $\tilde{\mathcal{G}}$ follows the configuration model. 

Consider the corresponding normalized Laplacian matrix $\tilde{M}$ for the pruned graph $\tilde{\mathcal{G}}$
$$\tilde{M}=\sqrt{\frac{\tilde{D}}{n}}\tilde{\Delta}^{-\frac{1}{2}}\left(\tilde{A}-\left[\frac{\tilde{d_i}\tilde{d_j}}{\tilde{D}-1}\right]_{ij}\right)\tilde{\Delta}^{-\frac{1}{2}}$$
where $\tilde{D}, \tilde{\Delta}, \tilde{A}, \tilde{d_i}$ are defined accordingly. 

Since the law for $\tilde{\mathcal{G}}$ follows the configuration model, we apply Theorem \ref{semicircle} and conclude that the empirical distribution of eigenvalues $\tilde{\mu}_n$ of $\tilde{M}$ weakly converges to the semicircle distribution in probability. In this proof, we will make connections between the normalized Laplacians $M$ and $\tilde{M}$ and show that their eigenvalue distributions are asymptotically the same. 

Note that $A$ is the adjacency matrix for $\mathcal{G}$ and $\tilde{A}$ is the adjacency matrix for $\tilde{\mathcal{G}}$. Therefore, $\tilde{A}$ is simply the resulting matrix by symmetrically setting at most $2\varepsilon n$ rows and columns in $A$ to be $0$. We use the index set $\mathcal{I}$ to denote the corresponding rows and columns of $0$'s. Here, we have $\tilde{A}^{(\mathcal{I})} = A^{(\mathcal{I})}$. Meanwhile, it is immediate that $\tilde{D} = (1 - o(1) ) D$.  Additionally, the matrix $\left[\frac{\tilde{d_i}\tilde{d_j}}{\tilde{D}-1}\right]_{ij}$ will not modify the eigenvalue distribution due to our observations in Remark \ref{rmk:rank1}. In other words, the Laplacians are connected in the following way,
\begin{align*}
    \left( \tilde{\Delta}^{\frac{1}{2}} \tilde{M} \tilde{\Delta}^{\frac{1}{2}} \right)^{(\mathcal{I})} &= (1-o(1)) \left( \Delta^{\frac{1}{2}} M \Delta^{\frac{1}{2}} \right)^{(\mathcal{I})}.
\end{align*}
Hence by applying $(\mathcal{I})$ to all matrices and moving terms,
\begin{align*}
    (\Delta^{(\mathcal{I})})^{-\frac{1}{2}} (\tilde{\Delta}^{(\mathcal{I})})^{\frac{1}{2}} \tilde{M}^{(\mathcal{I})} (\tilde{\Delta}^{(\mathcal{I})})^{\frac{1}{2}} (\Delta^{(\mathcal{I})})^{-\frac{1}{2}}  
    &= (1-o(1))  M^{(\mathcal{I})} .
\end{align*}
% \begin{align*}
%     \tilde{M}^{(\mathcal{I})} &= (1+o(1)) \tilde{\Delta}^{-\frac{1}{2}} \left( \Delta^{\frac{1}{2}} M \Delta^{\frac{1}{2}} \right)^{(\mathcal{I})} \tilde{\Delta}^{-\frac{1}{2}}\\
%     &= (1+o(1)) \tilde{\Delta}^{-\frac{1}{2}} (\Delta^{(\mathcal{I})})^{\frac{1}{2}} M^{(\mathcal{I})} (\Delta^{(\mathcal{I})})^{\frac{1}{2}}  \tilde{\Delta}^{-\frac{1}{2}}.
% \end{align*}
% By the Cauchy interlacing theorem, we know that the eigenvalues of $\tilde{A}$ interlace the eigenvalues of $A$. 
The most important part is how $\tilde{\Delta}^{(\mathcal{I})}$ differs from $\Delta^{(\mathcal{I})}$. We first establish an important observation. Denote the number of vertices whose degrees decrease by a factor of at least $\sqrt{\varepsilon}$ as 
$$K = \sum_{i=1}^n \mathbbm{1}_{\frac{d_i - \tilde{d}_i}{d_i} \geq \sqrt{\varepsilon}}.$$
Noting that the total number of deleted half-edges is less than $4\varepsilon C \sqrt{nD}$, we conclude the following
\begin{align*}
    K &\leq 3\varepsilon n + \sum_{i=1}^n \mathbbm{1}_{\frac{d_i - \tilde{d}_i}{d_i} \geq \sqrt{\varepsilon} \text{ and } d_i \geq C\sqrt{\frac{D}{n}}}\\
    &\leq 3 \varepsilon n + \sum_{i=1}^n \mathbbm{1}_{\frac{d_i - \tilde{d}_i}{C\sqrt{\frac{D}{n}}} \geq \sqrt{\varepsilon}}\\
    &\leq 3\varepsilon n + 4 \sqrt{\varepsilon} n
\end{align*}
where $3\varepsilon n$ accounts for the portion of vertices that are entirely removed. The third inequality holds because otherwise, the total number of removed half-edges will be greater than $4\sqrt{\varepsilon} n \times \sqrt{\varepsilon} C \sqrt{D/n} = 4 \varepsilon C \sqrt{nD}$, which forms a contradiction. Therefore, $O(\sqrt{\varepsilon }n)$ entries in the diagonal matrix $\tilde{\Delta}^{(\mathcal{I})}$ will be less than a factor of $1 - \sqrt{\varepsilon}$ of the corresponding entries in $\Delta^{(\mathcal{I})}$. The rest $n - O(\sqrt{\varepsilon }n)$ entries in $\tilde{\Delta}^{(\mathcal{I})}$ will be bounded above this cutoff. 

We use $\mathcal{J}$ to denote the index set of the $O(\sqrt{\varepsilon }n)$ rows and columns where their corresponding entries in $\tilde{\Delta}^{(\mathcal{I})}$ are respectively less than a factor of $1 - \sqrt{\varepsilon}$ of their corresponding entries in $\Delta^{(\mathcal{I})}$. We further take off this part of the rows and columns with respect to the equation above and get
\begin{equation}\label{transfer}
    (I - \Delta_\varepsilon)^{\frac{1}{2}} \tilde{M}^{(\mathcal{I},\mathcal{J})} (I - \Delta_\varepsilon)^{\frac{1}{2}}
    = (1-o(1))  M^{(\mathcal{I},\mathcal{J})} 
\end{equation}
where $\Delta_\varepsilon$ is a diagonal matrix where all entries are within the range of $(0,\sqrt{\varepsilon})$. We want to show that the multiplication of $(I - \Delta_\varepsilon)^{\frac{1}{2}}$ on both sides does not vary the eigenvalues much. 

Consider a symmetric matrix $B(t)$, its eigenvalue $\lambda (t)$, and its normalized eigenvector $v(t)$ which all depend on a variable $t$. By taking derivatives on both sides of the equation $B v = \lambda v$, we have that $\dot{B}v + B \dot{v} = \dot{\lambda} v + \lambda \dot{v}$. Taking dot product with $v$ on both sides and noting that $v$ is perpendicular to $\dot{v}$, we have $\dot{\lambda} = \langle v, \dot{B}v \rangle$. If we substitute $B(t) = (I - \Delta_\varepsilon t)^{\frac{1}{2}} \tilde{M}^{(\mathcal{I},\mathcal{J})} (I - \Delta_\varepsilon t )^{\frac{1}{2}}$, we know that 
\begin{align*}
    |\dot{\lambda}| &= |\langle v, -\Delta_\varepsilon (I - \Delta_\varepsilon t)^{-\frac{1}{2}} \tilde{M}^{(\mathcal{I},\mathcal{J})} (I - \Delta_\varepsilon t)^{\frac{1}{2}} v \rangle |\\
    &= |\langle v, - \Delta_\varepsilon (I - \Delta_\varepsilon t)^{-1} \lambda v \rangle |\\
    &\leq |\lambda | \cdot \norm{v}^2 \cdot \norm{\Delta_\varepsilon (I - \Delta_\varepsilon t)^{-1}} \\
    &\leq 2\sqrt{\varepsilon} |\lambda |.
\end{align*}
It immediately follows that
$$|\lambda (1) - \lambda (0)| \leq \int_0^1 |\dot{\lambda}| dt \leq 2\sqrt{\varepsilon} |\lambda| \cdot 1 \to 0.$$
Recall that $\lambda(1)$ corresponds to a eigenvalue for the matrix $B(1) = (I - \Delta_\varepsilon)^{\frac{1}{2}} \tilde{M}^{(\mathcal{I},\mathcal{J})} (I - \Delta_\varepsilon)^{\frac{1}{2}}$ and that $\lambda(0)$ corresponds to a eigenvalue for the matrix $B(0) = \tilde{M}^{(\mathcal{I},\mathcal{J})}$. 

Suppose we take an arbitrary Lipschitz function $f$. Note that the eigenvalue distribution of $\tilde{M}$ weakly converges to the semicircle distribution in probability. By Cauchy interlace theorem, the eigenvalue distribution of $\tilde{M}^{(\mathcal{I},\mathcal{J})}$ also weakly converges to the semicircle distribution in probability since we take off at most $|\mathcal{I}|+|\mathcal{J}| = O(\varepsilon n) + O(\sqrt{\varepsilon}n)$ entries. Therefore, we know that, with probability $1 - o(1)$,
$$\frac{1}{n} \sum_{i=1}^n f(\lambda_i(0)) \to \int f(x) \rho_{SC}(x) dx$$
Particularly, if we simply take the second moment, we have that 
$$\frac{1}{n} \sum_{i=1}^n \lambda_i(0)^2 \to \int x^2 \rho_{SC}(x) dx = C$$
where $C$ is the second moment of a semicircle distribution. Moreover, if we consider the difference between $\lambda (1)$ and $\lambda (0)$, with probability $1 - o(1)$,
\begin{align*}
    \frac{1}{n} \sum_{i=1}^n f(\lambda_i(1)) - \frac{1}{n} \sum_{i=1}^n f(\lambda_i(0))
    &=  \frac{1}{n} \sum_{i=1}^n f(\lambda_i(1)) - f(\lambda_i(0))\\
    &\leq \frac{1}{n} \norm{f'}_\infty \cdot \sum_{i=1}^n |\lambda_i(1) - \lambda_i(0)| \quad \text{ because $f$ is Lipschitz} \\
    &\leq \frac{1}{n} \norm{f'}_\infty \cdot \sum_{i=1}^n 2 \sqrt{\varepsilon} |\lambda_i (0)|\\
    &\leq  2 \sqrt{\varepsilon} \norm{f'}_\infty \cdot \sqrt{\frac{\sum_{i=1}^n |\lambda_i (0)|^2}{n}}\\
    &\to  2 \sqrt{\varepsilon} \norm{f'}_\infty \cdot \sqrt{C}\\
    &\to 0
\end{align*}
Therefore, we have that, with probability $1 - o(1)$,
\begin{align*}
    \frac{1}{n} \sum_{i=1}^n f(\lambda_i(1)) &= 
    \frac{1}{n} \sum_{i=1}^n f(\lambda_i(0)) + \frac{1}{n} \sum_{i=1}^n f(\lambda_i(1)) - \frac{1}{n} \sum_{i=1}^n f(\lambda_i(0))
    \to \int f(x) \rho_{SC}(x) dx,
\end{align*}
which implies that the eigenvalue distribution of $B(1) = (I - \Delta_\varepsilon)^{\frac{1}{2}} \tilde{M}^{(\mathcal{I},\mathcal{J})} (I - \Delta_\varepsilon)^{\frac{1}{2}}$ weakly converges to the semicircle distribution in probability.

Moreover, by \eqref{transfer}, the eigenvalue distribution of $M^{(\mathcal{I},\mathcal{J})}$ weakly converges to the semicircle distribution in probability. Again by Cauchy interlace theorem, the eigenvalue distribution of $M$ weakly converges to the semicircle distribution in probability, since we simply add at most $|\mathcal{I}| + |\mathcal{J}| = O(\varepsilon n) + O(\sqrt{\varepsilon}n)$ entries to $M^{(\mathcal{I},\mathcal{J})}$. Therefore, we know that the eigenvalue distribution of $M$ weakly converges to the semicircle distribution in probability.

\section{Computational Verification} \label{numerical}

\begin{figure}
    \begin{minipage}[c]{0.49\textwidth}
        \includegraphics[width=\linewidth]{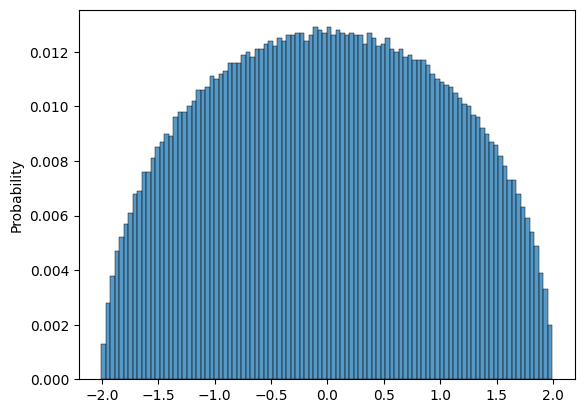}
        \caption{Spectral distribution for $d_1=100$, $d_2=500$, $n=10^5$.}
        \label{fig:semicircle}
    \end{minipage}
    \hfill
    \begin{minipage}[c]{0.49\textwidth}
        \includegraphics[width=\linewidth]{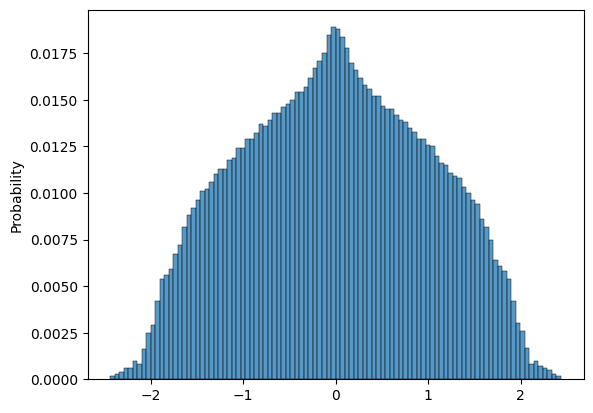}
        \caption{Spectral distribution
        for $d_1=10$, $d_2=200$, $n=10^5$.}
        \label{fig:notsemicircle}
    \end{minipage}%
\end{figure}

To computationally verify our results, we used Python along with numerical computing libraries NumPy and SciPy to produce distributions of eigenvalues\footnote{The code is available at https://github.com/kevinmli/Spectral-Distribution-Given-Degree-Sequences.}. In particular, we randomly assigned each of $n=10,000$ vertices one of two degrees (denoted $d_1,d_2$), then generated a random graph based on these fixed degree sequences. Figure \ref{fig:semicircle} and Figure \ref{fig:notsemicircle} show the spectral distributions for different values of $(d_1,d_2)$. In particular, Figure \ref{fig:semicircle} demonstrates the semicircle law when our assumption that $d_1,d_2>\sqrt{D/n}$ is satisfied while Figure \ref{fig:notsemicircle} demonstrates that when our assumption is not satisfied, the matrix can not follow the semicircle law.

In Figure \ref{fig:semicircle}, the fixed degree sequence matches our required conditions ($d_1,d_2>\sqrt{D/n}$), resulting in the semicircle distribution and supporting our conclusion. Note that Figure \ref{fig:notsemicircle} is not a semicircle distribution. Since we used $d_1=10$ and $d_2=200$, we have $D/n\approx 105$ so $\sqrt{D/n}\approx 10.25$. Since in this case $d_1 = 10 < \sqrt{D/n} \approx 10.25$, this provides evidence that our assumption that $d_i>\sqrt{D/n}$ is necessary.

%\section{Conclusion}
%In this paper, we investigate random graphs with given degree sequence $(d_1, d_2, \cdots , d_n)$ under mild conditions that $d_i = \omega (\sqrt{D/n})$ for all $i$. By applying the method of Fourier transform to deal with correlations and various transformations on the graphs constructed out of $k$-th moment of matrices, we proved the semicircle law for the normalized Laplacian matrix. Additionally, by using the birth and death chain to model the graph pruning process, we extend the assumption slightly to show the equivalent condition to the semicirclw law. 

\vspace{0.5in}
\noindent {\bf Acknowledgments} 
JH is supported by NSF grant DMS-2331096, and the Sloan research fellowship.

\bibliography{References.bib}
\bibliographystyle{abbrv}

\appendix
\section{Proof of \Cref{e:weak_converge}}
\label{sec:appendix1}

\begin{proof}[Proof of \Cref{e:weak_converge}]
    Combining \Cref{e:moment_converge} and Theorem \ref{concentration}, we have for any positive integer $k$
    $$ \lim_{n\rightarrow \infty} \p \biggl( \Big| \int x^k d\mu_n - \int x^k d\rho_{\rm sc} \Big| > \varepsilon \biggr) \to 0.
    $$
    Therefore, for any polynomial $p(x)$, we have 
    \begin{align}
        \label{e:poly_converge}
        \lim_{n\rightarrow \infty} \p \biggl( \Big| \int p(x) d\mu_n - \int p(x) d\rho_{\rm sc} \Big| > \varepsilon \biggr) \to 0.
    \end{align} 

    Now consider any continuous bounded function $f$ such that $\sup_{x\in \bR}|f(x)|\leq M$, and any arbitrarily small $\varepsilon > 0$. We want to show that 
    $$ \lim_{n\rightarrow \infty} \p \biggl( \Big| \int f(x) d\mu_n - \int f(x) d\rho_{\rm sc} \Big| > \varepsilon \biggr) \to 0.
    $$
    Fix an arbitrarily large number $K$. We define the truncated function $\tilde{f}$ as
    $$\tilde{f}(x) = \begin{cases} 
            f(x) & \text{if $x \in [-K,K]$}, \\
            0 & \text{otherwise}.
        \end{cases}$$
    Now observe that 
    $$\mu_n( [-K, K]^C ) \leq \frac{\int x^2 d\mu_n}{K^2} \leq \frac{2C_1}{K^2}$$
    where $C_1$ is Catalan number with parameter $1$. We take large $K$ such that $\frac{2C_1M}{K^2} < \frac{\varepsilon}{5}$. Therefore, we know 
    \begin{equation}
    \label{trun1}
        \Big| \int f d\mu_n - \int \tilde{f} d\mu_n \Big| = \Big| \int_{[-K,K]^C} f d\mu_n \Big| \leq \sup_{x\in [-K,K]^C} |f(x)| \cdot \mu_n ([-K, K]^C) \leq M \frac{2C_1}{K^2} < \frac{\varepsilon}{5}.
    \end{equation}
    By Stone-Weierstrass theorem, on the finite interval of $[-K,K]$, we can approximate any continuous function by polynomials to arbitrary precision. Therefore, we may choose a polynomial $p(x)$ which controls the difference between $\tilde{f}$ and $p$ on the interval of $[-K,K]$. On the interval of $[-K,K]^C$, a similar argument as in \eqref{trun1} can be applied. Therefore, by combining two arguments and choosing a polynomial with finer precision, we have
    \begin{equation}
    \label{trun2}
         \Big| \int \tilde{f} d\mu_n - \int p d\mu_n \Big| < \frac{\varepsilon}{5}  
    \end{equation}
    and 
    \begin{equation}
    \label{trun3}
     \Big|  \int p d\rho_{\rm sc} - \int \tilde{f} d\rho_{\rm sc}\Big| < \frac{\varepsilon}{5}. 
    \end{equation}
 We notice the fact that
    \begin{equation}
    \label{trun5}
    \Big| \int \tilde{f} d\rho_{\rm sc} - \int f d\rho_{\rm sc} \Big| = 0
    \end{equation}
    because $\rho_{\rm sc}$ has a support of $[-2,2]$. Finally, by observations \eqref{e:poly_converge}, \eqref{trun1}, \eqref{trun2}, \eqref{trun3}, and \eqref{trun5}, we conclude that
    \begin{align*}
        & \lim_{n\rightarrow \infty} \p \biggl( \Big| \int f(x) d\mu_n - \int f(x) d\rho_{\rm sc} \Big| > \varepsilon \biggr) \\
        = &\lim_{n\rightarrow \infty} \p \biggl( 
        \Big| \int fd\mu_n - \int \tilde{f}d\mu_n +  \int \tilde{f}d\mu_n - \int pd\mu_n + \int pd\mu_n - \int pd\rho_{\rm sc} 
        + \int pd\rho_{\rm sc} - \int \tilde{f}d\rho_{\rm sc} \\
        &+ \int \tilde{f}d\rho_{\rm sc} - \int fd\rho_{\rm sc} \Big| > \varepsilon
        \biggr)\\
        \leq &\lim_{n\rightarrow \infty} \p \biggl( \Big| \int f d\mu_n - \int \tilde{f} d\mu_n \Big| > \frac{\varepsilon}{5}\biggr) +  \lim_{n\rightarrow \infty} \p \biggl( \Big| \int \tilde{f} d\mu_n - \int p d\mu_n \Big| > \frac{\varepsilon}{5}\biggr)\\ &+ \lim_{n\rightarrow \infty} \p \biggl( \Big| \int p d\mu_n - \int p d\rho_{\rm sc} \Big| > \frac{\varepsilon}{5}\biggr) 
        + \lim_{n\rightarrow \infty} \p \biggl( \Big|  \int p d\rho_{\rm sc} - \int \tilde{f} d\rho_{\rm sc} \Big| > \frac{\varepsilon}{5}\biggr) \\
        &+ 
        \lim_{n\rightarrow \infty} \p \biggl( \Big| \int \tilde{f} d\rho_{\rm sc} - \int f d\rho_{\rm sc} \Big| > \frac{\varepsilon}{5}\biggr) \to 0
    \end{align*}
    which concludes our proof.
\end{proof}

\end{document}